\documentclass[preprint,12pt]{elsarticle}

\usepackage{amssymb}
\usepackage{amsthm}

\usepackage{amsmath} 
\DeclareMathOperator{\sign}{sign}
\usepackage{tikz}
\usetikzlibrary{arrows, arrows.meta} 
\usetikzlibrary{decorations.pathreplacing}
\usetikzlibrary{decorations.pathmorphing,patterns}
\usetikzlibrary{calc,patterns,decorations.markings}
\usetikzlibrary{positioning}
\usepackage{graphicx}
\usepackage{caption}
\usepackage{changepage}
\usepackage{adjustbox}
\usepackage{lipsum}
\usepackage{subcaption}

\usepackage{calrsfs}
\usepackage[showframe=false]{geometry}
\DeclareMathAlphabet{\pazocal}{OMS}{zplm}{m}{n}

\usepackage[section]{placeins}
\newenvironment{rcases}
  {\left.\begin{aligned}}
  {\end{aligned}\right\rbrace}

\theoremstyle{definition}

\theoremstyle{remark}

\graphicspath{
{Figs/}}

\usepackage{float}
\begin{document}
\begin{frontmatter}



\title{Novel, simple and robust contact-discontinuity capturing schemes for high speed compressible flows}


\author[1,2]{Ramesh Kolluru\corref{cor1}}
\ead{kollurur@iisc.ac.in}
\ead{rameshsurya335@gmail.com}
\cortext[cor1]{Corresponding author} 
\author[1,4]{N. Venkata Raghavendra}
\ead{venkata.r.nandagiri@gmail.com}
\author[1]{S. V. Raghurama Rao}
\ead{raghu@iisc.ac.in}
\ead{svraghuramarao@gmail.com}
\author[3]{G.N.Sekhar}
\ead{drgns.maths@bmsce.ac.in}
\ead{gnsnms@gmail.com}

\address[1]{Department of Aerospace Engineering, Indian Institute of Science, Bangalore\vspace{0.1cm}}
\address[2]{Department of Mechanical Engineering, B.M.S.College of Engineering, Bangalore\vspace{0.1cm}}
\address[3]{Department of Mathematics, B.M.S.College of Engineering, Bangalore\vspace{0.1cm}}
\address[4]{Department of Aeronautical Engineering, Anna Saheb Dange College of Engineering, Sangli, Maharastra}

\begin{abstract}
The nonlinear convection terms in the governing equations of compressible fluid flows are hyperbolic in nature and are nontrivial for modelling and numerical simulation. Many numerical methods have been developed in the last few decades for this purpose and are typically based on Riemann solvers, which are strongly dependent on the underlying eigen-structure of the governing equations. Objective of the present work is to develop simple algorithms which are not dependent on the eigen-structure and yet can tackle easily the hyperbolic parts. Central schemes with smart diffusion mechanisms are apt for this purpose.  For fixing the numerical diffusion, the basic ideas of satisfying the Rankine-Hugoniot (RH) conditions along with generalized Riemann invariants are proposed. Two such interesting algorithms are presented, which capture grid-aligned steady contact discontinuities exactly and yet have sufficient numerical diffusion to avoid numerical shock instabilities. Both the algorithms presented are robust in avoiding shock instabilities, apart from being accurate in capturing contact discontinuities, do not need wave speed corrections and are independent of eigen-strutures of the underlying hyperbolic parts of the systems.  
\end{abstract} 
\begin{keyword}
Compressible flow solvers, Rankine-Hugoniot jump condition, Riemann invariants, robust central schemes, accurate contact discontinuity capturing, eigen-structure independence 



\end{keyword}

\end{frontmatter}


\section{Introduction}
Development of numerical algorithms for simulating compressible fluid flows is an active area of research. The quest to develop simple, robust and low numerical diffusion algorithms has been a continuing feature of research in CFD in the past several decades. For a detailed review of these schemes the reader is referred to \cite{hll,jameson3,jameson4,maccormack1,roe1,roe2,leer1,Jaisankar_SVRRao,Jaisankar_SVRRao1,Venkat_thesis,TORO_1,Laney,Leveque}. The following features are worth considering while developing new schemes for hyperbolic systems representing gas dynamics.  
 \begin{itemize}
 \item Exact capturing of steady discontinuities.
 \item Minimum numerical diffusion without violating entropy conditions. 
 \item Avoiding shock instabilities.
 \item Eigen-structure independency.
 \item Simplicity of the algorithm.
 \end{itemize}
Most popular algorithms depend on Riemann solvers and eigen-structure.  Though some of them can capture grid aligned steady shocks or contact-discontinuities exactly, they often produce unphysical phenomena like carbuncle shocks, kinked Mach stems, odd-even decoupling, and violation of entropy conditions, because of inherent low numerical diffusion present in them. Researchers in the recent past have focussed on algorithms which can avoid these anomalies and the search for an ideal scheme is still continuing. In this quest, it will be advantageous to incorporate the physical and mathematical features characteristic of the nonlinear propagating waves. Out of the three nonlinear waves, the shock waves must satisfy the Rankine-Hugoniot conditions, the rarefaction waves must satisfy the Generalized Riemann Invariants (GRI) and the contact-discontinuities must satisfy both \cite{TORO_1}. In the present work we use both the above criteria to develop two new, simple and robust algorithms for Euler and Navier-Stokes equations.  
The rest of the paper is organized as follows. In section \ref{sec:geqns} a brief introduction for the governing equations and time discretization in finite volume framework is presented.  The description of the new algorithms is presented in sections \ref{sec:RICCA} and section \ref{sec:movers+}.  In Sections \ref{sec:results} and \ref{sec:viscous_results}, the results produced by these new robust algorithms for various 1D and 2D bench-mark test cases both for Euler and Navier Stokes equations are presented, followed by a summary in the last section.    

\section{Governing equations} \label{sec:geqns}
\begin{figure}[h!]
\begin{center}
\begin{tikzpicture}[scale = 1.3]
\draw (0,1)node (xaxis) [below] {(j-1)}  -- (2,1)node (xaxis) [below] {(j)} -- (4,1)  node (xaxis) [below] {(j+1)};
\draw [red](2,0.5) node (xaxis) [below] {$L$}; 
\draw [red](4,0.5) node (xaxis) [below] {$R$};
\draw (2,1) node (xaxis)[above] {$F_L = F(U_L)$}; 
\draw (4,1) node (xaxis)[above] {$F_R = F(U_R)$};
\draw [red,thick,dashed](1,0) node (yaxis) [below] {$j-\frac{1}{2}$}-- (1,2);
\draw [red,thick,dashed](3,0) node (yaxis) [below] {$j+\frac{1}{2}$}-- (3,2);
\draw[red](3,2) node (yaxis) [above] {$F_\mathrm{I}$};
\draw [blue,fill] ( 1,1) circle [radius=0.03] ;
\draw [blue,fill] ( 3,1) circle [radius=0.03] ;
\draw [blue,fill] ( 0,1) circle [radius=0.03] ;
\draw [blue,fill] ( 2,1) circle [radius=0.03] ;
\draw [blue,fill] ( 4,1) circle [radius=0.03] ;
\end{tikzpicture}    
\caption{Typical finite volume in 1D}
\label{Interfaceflux1}
\end{center}
\end{figure}
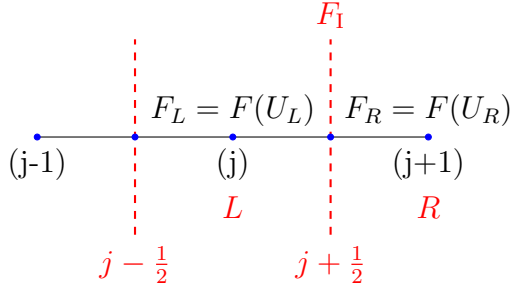

The basic equations which govern the dynamics of fluid flow are derived from conservation laws of mass, momentum and energy. The integral form of non-dimensional Navier-Stokes equations, with Fourier law of heat conduction, Sutherland's law and the expression for the viscous stresses are as given in (\ref{NonDim_NSeq}).
\begin{equation}
  \begin{aligned} \label{NonDim_NSeq}
\frac{d}{dt}\int_{\Omega} \rho ~d\Omega + \int_{S} \rho  V_n~dS &= 0, \\ 
V_n &= \vec{V}\cdot\hat n,\\ 
\frac{d}{dt}\int_{\Omega} \rho\vec{V} ~d\Omega + \int_{S} \left( (\rho \vec{V}) V_n  - p \hat n\right)~ dS  &= \frac{1}{Re}\int_{s}  \bf{\tau} \cdot d S, \\
\tau &=  \mu \left( {\frac{\partial u_i}{\partial x_j} + \frac{\partial u_j}{\partial x_i} - \frac{2}{3}\delta_{ij}\frac{\partial u_k}{\partial x_k}}\right), \\
\frac{d}{dt}\int_{\Omega} \rho E_{t} ~d\Omega +\int_{S} (\rho E_{t} + p)V_n~ dS &= \int_{S} Q d S + \int_{\Omega} \nabla \cdot {\left(\bf{\tau} \cdot \vec{v}\right) d \Omega},\\
Q = \frac{-\mu}{(\gamma -1) Re M^2_{ref} Pr } \frac{\partial T}{\partial x_i},
p = \frac{\rho T}{\gamma M^2_{ref}},
  \mu &= T^{\frac{3}{2}}\left(\frac{1 + \frac{T_s}{T_{ref}}}{T + \frac{T_s}{T_{ref}}}\right).
  \end{aligned}
  \end{equation}
In the above equations $\rho$ is the density of the fluid, $\vec{V}$ is velocity vector, $p$ is the thermodynamic pressure, and $E_t$ refers to the total energy per unit mass defined as $E_t = c_v T + \frac{|\vec{V}|^2}{2}= \frac{p}{\rho\left(\gamma - 1\right)} + \frac{|\vec{V}|^2}{2}$, with perfect gas EOS given by
  \begin{align} 
p &= \rho R T \label{Perfect_gas_EOS},
  \end{align}  
Neglecting viscous terms in (\ref{NonDim_NSeq}) results in the governing equations for inviscid compressible fluid flows, known as  Euler equations.  
These system equations are nonlinear in nature and therefore analytical solutions are hard to obtain. The nonlinearity of the convection terms leads to nonlinear waves (compression waves and expansion waves) and may eventually lead to discontinuities (shock waves and contact-discontinuities) even if the initial conditions are smooth. Once the discontinuities appear, the differential form of equations is no longer valid and the description needs to be shifted to weak form of conservation laws, which also yields the integral form.  The \emph{weak solution} refers to piecewise smooth solutions with discontinuities in between.  In this work cell centred finite volume framework is used to discretize both  Euler equations and NS equations. The basic governing equations described in (\ref{NonDim_NSeq}) can thus be rewritten in compact notation as 
\begin{subequations}
\begin{align}
\frac{d \overline{U}}{dt} = - R, \textrm{where} \
\label{Residue}
R &= \frac{1}{\Omega} \left[\sum_{i=1}^{N}{F_c \cdot \hat{n} ~dS} - \sum_{i=1}^{N}{F_v \cdot \hat{n}~dS}\right],\\
\label{averageU}
 \overline{U} &= \frac{1}{\Omega}\int_{\Omega} U d\Omega~.
\end{align}
\end{subequations}where $U$ is conserved variable vector (the bar representing a cell-integral average), $F_c$ is convective flux vector and $F_v$ is viscous flux vector  on an interface, $R$ represents  net flux from a given control volume, $\Omega$ the cell volume and $N$ represents the number of control surfaces for a given control volume.
The convective flux on any interface of a control volume as shown in figure (\ref{Interfaceflux1}) for any stable scheme can be written as the sum of an average flux across the interface and a numerical dissipative flux as given in \eqref{ConvectiveDiffusiveFlux}.
\begin{equation}
\centering
\label{ConvectiveDiffusiveFlux}
 F_\mathrm{I} = \frac{1}{2}\left[ F_L + F_R\right] - d_\mathrm{I} ;~~
d_\mathrm{I}  = \frac{ \alpha_{\mathrm{I}}}{2}\left(U_R - U_L\right)
\end{equation}
where $\alpha_{\mathrm{I}}$ is the  coefficient of numerical diffusion and $d_\mathrm{I}$ represents the numerical dissipative flux. In the present work the coefficient of numerical diffusion is fixed by the new algorithms RICCA and MOVERS+ which are explained in sections (\ref{sec:RICCA} $\&$ \ref{sec:movers+}).  
The discretization of viscous flux $F_v$ requires an evaluation of the second gradient of velocities and gradient of temperature at the centroid of the control volume.  In the present work Green-Gauss (GG) based method is used to evaluate the gradients on a structured grid, using a diamond structure in a co-volume.  The details of the viscous flux discretization are given in \cite{Ramesh_Thesis}.   
Once the spatial discretization is done (as described in the following subsection), the PDEs get converted to ODEs and then time discretization can be carried out using popular methods for solving ODEs such as Euler method or Runge Kutta Methods.  For a particular control volume, the first order (1O) Euler method is given by (\ref{1OEULER}) as 
\begin{equation}
\begin{aligned}\label{1OEULER}
\frac{\Delta \overline{U}_{ij}}{\Delta t} = -R^n_{ij}, \\ 
\overline{U}_{ij}^{n+1} = \overline{U}_{ij}^n - \Delta t R^n_{ij}.
\end{aligned}
\end{equation}
Higher order time accuracy is acheived using a third order (3O) Runge-Kutta method, as given in (\ref{RK3}).
\begin{equation}
\begin{aligned}\label{RK3}
\overline{U}^1_{ij} = \overline{U}^n_{ij} - \Delta t R^n_{ij}\left(\overline{U}^n_{ij}\right) \\ 
\overline{U}^2_{ij} = \frac{1}{4}\overline{U}^1_{ij}  + \frac{3}{4} \overline{U}^n_{ij} - \frac{1}{4}\Delta t R_{ij}\left(\overline{U}^1_{ij}\right) \\ 
\overline{U}^{n+1}_{ij} = \frac{2}{3}\overline{U}^2_{ij}  + \frac{1}{3}\overline{U}^n_{ij} - \frac{2}{3}\Delta t R_{ij}\left(\overline{U}^2_{ij}\right) 
\end{aligned}
\end{equation} 

\section{Riemann Invariant based Contact-discontinuity Capturing Algorithm (RICCA)}\label{sec:RICCA}
In this section a novel scheme is presented in which the effect of  Generalized Riemann Invariants is utilized in the discretization process, leading to a scheme which captures steady contact discontinuities exactly.  
 \subsection{Generalised Riemann Invariants (GRI)} \label{GRI}
 The concept of GRI is briefly introduced here (for a more detailed explanation see \cite{Jeffrey,TORO_1,Sonar}).  Consider a general quasi-linear hyperbolic system as given by  (\ref{qlcf}).  
 \begin{align}
 \label{qlcf}
\frac{\partial U}{\partial t} + A(U)\frac{\partial U}{\partial x}=0,\\ 
U = [U_1,U_2,\cdots ,U_m]^T,
\end{align}
where $U$ represents the conserved variable vector of the hyperbolic system. Of the $m$ waves associated with the system (\ref{qlcf}) for the $i^{th}$ characteristic field associated with eigenvalue $\lambda_i$, corresponding right eigenvector is given by 
\begin{equation}
\label{ieigenvector}
\bf{R^i} = \left[ \bf{r_1^i, r_2^i \cdots r_m^i} \right]^T
\end{equation}
The Generalised Riemann Invariants are relations that hold true across expansion waves and contact-discontinuites. This can be mathematically written as 
\begin{equation}
\label{gri}
\frac{dU_1}{\bf{r_1^i}} = \frac{dU_2}{\bf{r_2^i}}= \cdots = \frac{dU_m}{\bf{r_m^i}}
\end{equation}
These equations relate ratios of $dU_j$ to the respective component $r_i^j$ of the right eigenvector $R_i^j$, corresponding to an eigenvalue $\lambda_{i}$.  Here, the above relations (GRIs) are utilized in developing a new algorithm which can recognize contact-discontinuities and the algorithm is expected to be accurate enough for flow simulations. The above ideas will be incoporated in a simple central discretization framework in the finite volume method, avoiding Riemann solvers, field-by-field decompositions and complicated flux splittings. This is achieved by fixing the coefficient of numerical diffusion in a generic expression for the interface flux based on the above criteria.
\subsection{A central solver based on GRI} 
 The finite volume update formula for Euler equations is given by \eqref{1DEulerupdate} with interface flux evaluated as in  \eqref{interfaceflux1}.
\begin{align}
\label{1DEulerupdate}
 \bar{U}^{n+1}_j &=    \bar{U}^{n}_j - \frac{\Delta t}{\Delta x}\left[ F^n_{j + \frac{1}{2}} - F^n_{j - \frac{1}{2}} \right] \\
\label{interfaceflux1}
 F_{j \pm \frac{1}{2}} = F_I\left(U_L,U_R \right) &= \frac{1}{2}\left[ F(U_L) + F(U_R) \right] - \Delta F_{num} 
\end{align}
where the first term on the right hand side is an average flux from the left (L) and the right (R) states and $\Delta F_{num}$ is a flux difference  representing numerical diffusion. This numerical diffusion is modeled as follows. 
\begin{equation}  
\Delta F_{num} = \left( \displaystyle \frac{\Delta F}{\Delta U} \right)_{num} \Delta U = \alpha_{num} \Delta U = \alpha_{I} \Delta U 
\end{equation} 

\begin{figure}[h!]
\centering
\includegraphics[trim = 4mm 4mm 10mm 10mm,width=0.28\linewidth]{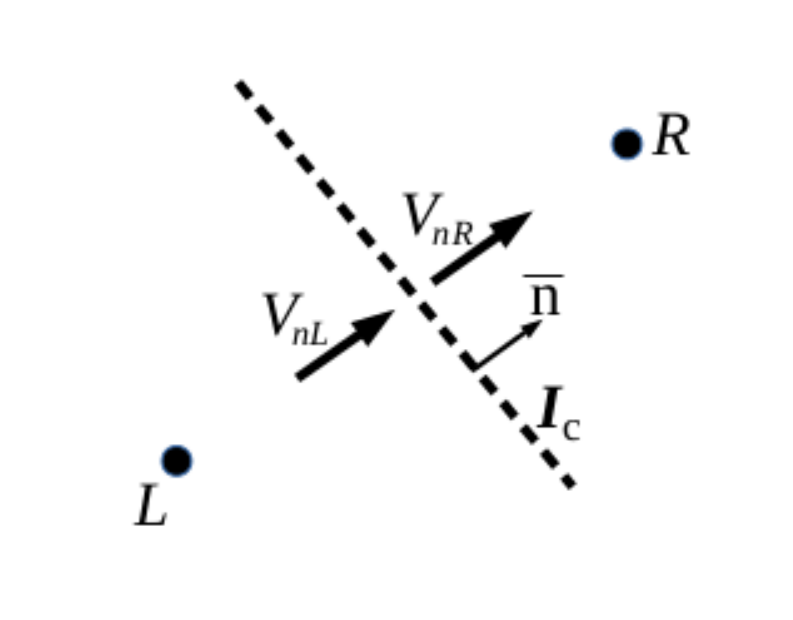}
\caption{Locally one-dimensional flow at cell interface $I_c$}
\label{fig:locally_1d_flow}
\end{figure}

The coefficient of numerical diffusion, $\alpha_{num}$, is modeled using a diagonal matrix such that
\begin{eqnarray} 
 \alpha_{I} = \alpha_{num} = \left[\begin{array}{ccc} \alpha_{1} & 0 & 0 \\ 0 & \alpha_{2} & 0 \\ 0 & 0 & \alpha_{3} \end{array} \right] = 
 \left[ \begin{array}{ccc} \left| \displaystyle \frac{\Delta F_{1}}{\Delta U_{1}}\right| & 0 & 0 \\ 
 0 & \left| \displaystyle \frac{\Delta F_{2}}{\Delta U_{2}} \right| & 0 \\ 
 0 & 0 & \left| \displaystyle \frac{\Delta F_{3}}{\Delta U_{3}} \right|  
 \end{array} \right] 
\end{eqnarray} 
  Various numerical schemes differ in the way this wave speed or the coefficient of numerical diffusion is determined. The basic idea of the present work is to use GRI across the interface to determine the coefficient of diffusion, $\alpha_{I}$. As shown  \cite{TORO_1,Sonar}, a contact-discontinuity separates two states in a linearly degenerate field and constancy of the GRI holds good across this wave, apart from R-H conditions. The eigenvalue corresponding to the contact-discontinuity is given by $\lambda = u$, with the corresponding right eigenvector as in  \eqref{eigenvector_u}. 
\begin{align}
\label{eigenvector_u}
\mathbf{R}  =\left [ \begin{array}{c}
r_1 \\
r_2\\
r_3
\end{array} \right ]= \left [ \begin{array}{c}
1 \\
u \\
\frac{u^2}{2}
\end{array} \right ].
\end{align}
The GRIs (\ref{gri}) applied to the contact-discontinuity leads to the following ODEs:
\begin{equation}
\label{odes}
\frac{d \rho}{1} = \frac{d{(\rho u)}}{u} = \frac{d (\rho E)}{\frac{u^2}{2}}.
\end{equation}
Solving the above ODEs \eqref{odes} results in pressure and velocity being constant across a contact-discontinuity \cite{TORO_1,Sonar}.   Therefore the conditions for pressure and velocity across the cell interface can be written as \eqref{1d_contact_conditions}
\begin{equation}
\begin{aligned}\label{1d_contact_conditions}
\Delta{u}&=0 \ or \ u_{j}=u_{j+1} = u_{\raisebox{-2pt} {\scriptsize \emph{$I$}}}, \\
\Delta{p}&=0 \ or \ p_{j}=p_{j+1} = p_{\raisebox{-2pt} {\scriptsize \emph{$I$}}}.
\end{aligned}
\end{equation} 
 Using the conditions \eqref{1d_contact_conditions} in the expressions of $(\alpha_{I})_{l}, \quad l = 1,2,3$, one can obtain the following. \\
$l=1$: .
\begin{equation}
({\alpha_I})_1 = \bigg|\frac{(F_1)_{j+1}-(F_1)_{j}}{(U_1)_{j+1}-(U_1)_{j}}\bigg| = \bigg|\frac{\rho_{j+1}u_{j+1}-\rho_{j}u_{j}}{\rho_{j+1}-\rho_{j}}\bigg| = |u_{\raisebox{-2pt} {\scriptsize \emph{$I$}}}|
\end{equation}
$l=2$:
\begin{equation}
({\alpha_I})_2 = \bigg|\frac{(F_2)_{j+1}-(F_2)_{j}}{(U_2)_{j+1}-(U_2)_{j}}\bigg| = \bigg|\frac{p_{j+1}+\rho_{j+1}u^2_{j+1}-p_{j}-\rho_{j}u^2_{j}}{\rho_{j+1}u_{j+1}-\rho_{j}u_{j}}\bigg| = |u_{\raisebox{-2pt} {\scriptsize \emph{$I$}}}|
\end{equation}
$l=3$:
\begin{equation}
({\alpha_I})_3 = \bigg|\frac{(F_3)_{j+1}-(F_3)_{j}}{(U_3)_{j+1}-(U_3)_{j}}\bigg| = \bigg|\frac{p_{j+1}u_{j+1}+\rho_{j+1}u_{j+1}E_{j+1}-p_{j}u_{j}-\rho_{j}u_{j}E_{j}}{\rho_{j+1}E_{j+1}-\rho_{j}E_{j}}\bigg| = |u_{\raisebox{-2pt} {\scriptsize \emph{$I$}}}|
\end{equation}
From the above three expressions the coefficient of numerical diffusion determined to accurately capture contact-discontinuity located at the cell interface ${I} \equiv j+\frac{1}{2}$ is 
\begin{equation}
({\alpha}_{\mathrm{I}})_l = |u_{\raisebox{-2pt} {\scriptsize \emph{$I$}}}|, \quad l=1,2,3
\end{equation} 
leading a scalar numerical diffusion.  
Using (\ref{1d_contact_conditions}), the above coefficient of diffusion can be expressed in four different ways as
\begin{equation}
\label{1dcontactdiffusion}
({\alpha}_{\mathrm{I}})_l = |u_j| = |u_{j+1}| = \frac{|u_{j}|+|u_{j+1}|}{2} = max(|u_j|, |u_{j+1}|), \quad l=1,2,3.
\end{equation}
In a general multi-dimensional flow case, if a locally 1D flow is assumed at the cell interface as depicted in figure (\ref{fig:locally_1d_flow}),
\vspace{3mm}
the coefficient of numerical diffusion to accurately capture a contact-discontinuity can be expressed as 
\begin{equation}\label{general_contact_capture}
{\alpha}_{\raisebox{-2pt} {\scriptsize \emph{$I$}}} = |V_{nL}| = |V_{nR}| = \frac{|V_{nL}|+|V_{nR}|}{2} = max(|V_{nL}|, |V_{nR}|).
\end{equation}
Numerical experimentation has revealed that, this numerical diffusion evaluated by \eqref{1dcontactdiffusion} or \eqref{general_contact_capture}, though adequate in capturing the contact-discontinuities exactly, is not sufficient enough for the case of shocks being located at the cell interface. So in order to generalize the diffusion for any case the Riemann Invariant based Contact-discontinuity Capturing Algorithm (RICCA) is designed with the following coefficient of numerical diffusion:
\begin{equation}\label{accudisks_eval_alpha_euler_2d}
{\alpha}_ {\emph{$I$}} = \begin{cases}
                            \qquad \qquad \frac{|V_{nL}|+|V_{nR}|}{2}, \qquad \qquad \qquad \quad \text{if } |\Delta\mathbf{F}|<\delta \ \text{and} \ |\Delta\mathbf{U}| <\delta\\
                            \quad \qquad \qquad \qquad \qquad \qquad \qquad \qquad \qquad  \\
                            max(|V_{nL}|, |V_{nR}|) + sign(|\Delta p_{\raisebox{-2pt} {\scriptsize {$\mathrm{I}$}}}|) a_{\mathrm{I}} , \qquad \text{otherwise} \\
                            \qquad \qquad \qquad \qquad \quad
                           \end{cases}
\end{equation}
where $\delta$ is a small number and  $a_{\mathrm{I}} = \sqrt{\frac{\gamma p_{\raisebox{-2pt} {\scriptsize \emph{$I$}}}}{\rho_{\raisebox{-2pt} {\scriptsize \emph{$I$}}}}}$ is the speed of sound evaluated with the values at the interface given by
\begin{align}
\rho_{\raisebox{-2pt} {\scriptsize \emph{$I$}}} &= \frac{\rho_L+\rho_R}{2}, \\
p_{\raisebox{-2pt} {\scriptsize \emph{$I$}}} &= \frac{p_L+p_R}{2}, \\
\Delta{p}_{\raisebox{-2pt} {\scriptsize \emph{$I$}}} &= (p_R-p_L).
\end{align}
From (\ref{accudisks_eval_alpha_euler_2d}) it can be seen that for the case of a steady contact-discontinuity at the interface where $sign(|\Delta p_{\raisebox{-2pt} {\scriptsize \emph{$I$}}}|) = 0$, the coefficient ${\alpha}_{\raisebox{-2pt} {\scriptsize \emph{$I$}}}$ becomes identical to the expression in (\ref{general_contact_capture}) resulting in exact capturing of the steady contact-discontinuity. On the other hand, if a shock is located at the interface in which case $sign(|\Delta p_{\raisebox{-2pt} {\scriptsize \emph{$I$}}}|) = 1$, the expression for the coefficient in (\ref{accudisks_eval_alpha_euler_2d}) becomes $max(|V_{nL}|, |V_{nR}|) + a_{\mathrm{I}}$ which is a Rusanov (LLF) type diffusion and should be adequate near shocks. Even in the case of an expansion region with a sonic point ($M=1$) at the interface, $sign(|\Delta p_{\raisebox{-2pt} {\scriptsize \emph{$I$}}}|) = 1$, the expression for the coefficient in (\ref{accudisks_eval_alpha_euler_2d}) becomes $max(|V_{nL}|, |V_{nR}|) + a_{\mathrm{I}}$ which results in non-zero diffusion ensuring no expansion shocks. So, entropy violation is unlikely to occur. This design of the coefficient of numerical diffusion therefore does not require any entropy fix.

On the whole, the new central scheme RICCA:
\begin{itemize}
 \item can capture steady grid-aligned contact-discontinuities exactly,
 \item has sufficient numerical diffusion near shocks so as to avoid shock instabilities, and
 \item does not need entropy fix for at sonic points. 
 \item is not tied down to the eigen-structure and hence can be easily extended to any general equation of state, without modification. 
\end{itemize}
A similar strategy was introduced by N.Venkata Raghavendra in \cite{Venkat_thesis,Venkat_Arxiv} to design an accurate contact-discontinuity capturing discrete velocity Boltzmann scheme for inviscid compressible flows.  

\section{New central scheme, MOVERS+} 
The second of the two new algorithms presented in this paper is based on  subtantial modification of a central Rankine-Hugoniot solver developed by Jaisankar \& Raghurama Rao \cite{Jaisankar_SVRRao}, called as MOVERS (Method of Optimal Viscosity for Enhanced Resolution of Shocks).  This is first briefly reviewed in the following subsection, before introducing the new scheme, named as MOVERS+.  
\subsection{MOVERS} 
MOVERS \cite{Jaisankar_SVRRao} is a central scheme which can capture grid aligned steady shocks and contact discontinuities exactly, without numerical diffusion.  As it is is a central scheme, it avoids all the complications of Riemann solvers and is not tied to the eigen-structure of the underlying hyperbolic systems.  The accurate discontinuity capturing is achieved by enforcing the Rankine-Hugoniot jump condition directly in the discretization process. The basic idea of this algorithm is briefly explained in the following. 

Consider the Rankine-Hugoniot conditions, given by \eqref{RHCondition}
\begin{equation}
\label{RHCondition}
 \Delta F = s \Delta U, \quad \Delta(\cdot) = (\cdot)_R - (\cdot)_L
\end{equation}
where $s$ is the speed of the discontinuity, $F$ is flux vector and $U$ is the conserved variable vector.
\begin{figure}[h!]
\begin{center}
\begin{tikzpicture}[scale = 1.2]
\draw (0,1)node (xaxis) [below] {(j-1)}  -- (2,1)node (xaxis) [below] {(j)} -- (4,1)  node (xaxis) [below] {(j+1)};
\draw [-,red,thick] (2.0,1.0) node (xaxis) [above] {$L$}; 
\draw [-,red,thick] (4.0,1.0) node (xaxis) [above] {$R$}; 
\draw (2.7,1.7) node (xaxis)[above] {$s_I^- $};
\draw [->,red,thick](2.7,1.7)--(2.5,1.7) ;
\draw (3.3,1.7) node (xaxis)[above] {$s_I^+$};
\draw [->,red,thick](3.2,1.7)--(3.4,1.7) ;
\draw [red,thick](1,0) node (yaxis) [below] {$j-\frac{1}{2}$}-- (1,2);
\draw [blue,thick,dashed](2.9,0) -- (2.9,2);
\draw [red,thick](3,0) node (yaxis) [below] {$j+\frac{1}{2}$}-- (3,2);
\draw [blue,thick,dashed](3.1,0)  -- (3.1,2);
\draw[red,thick](3,2.1) node (yaxis) [above] {$F_\mathrm{I}$};
\draw [blue,fill] ( 1,1) circle [radius=0.03] ;
\draw [blue,fill] ( 3,1) circle [radius=0.03] ;
\draw [blue,fill] ( 0,1) circle [radius=0.03] ;
\draw [blue,fill] ( 2,1) circle [radius=0.03] ;
\draw [blue,fill] ( 4,1) circle [radius=0.03] ;
\end{tikzpicture}    
\caption{Shock located at interface}
\label{splitshockspeed}
\end{center}
\end{figure}
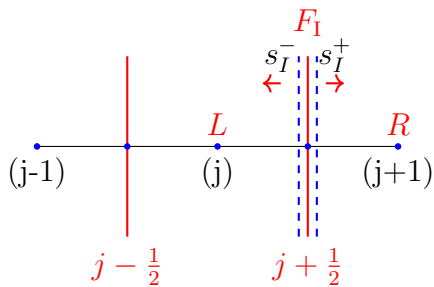
First, the speed of the discontinuity is split into a positive part (corresponding to a right-moving discontinuity) and a negative part (corresponding to a left-moving discontinuity), as shown in figure (\ref{splitshockspeed}) and as given in \eqref{splitwave}
\begin{align}
\label{splitwave}
s_I &= s_I^+ + s_I^- \\
s_I^{\pm} &= \frac{s_I \pm \left| s_I  \right|}{2}.
\end{align}
Using the above wave speed splitting, the RH condition at the interface is split into two parts as 
\begin{align}
\begin{rcases}
\label{SplitRHCondition}
F_{R} - F_{\mathrm{I}} = s_I^+ \Delta U, \\
F_{\mathrm{I}} - F_{L} = s_I^- \Delta U
\end{rcases}
\end{align}
These split RH conditions \eqref{SplitRHCondition} lead to the cell-interface flux as 
\begin{align}
\label{Final_Interfaceflux}
F_\mathrm{I} = \frac{F_{L} + F_{R}}{2} - \frac{\left| s_I \right|}{2} \Delta U 
\end{align}
Comparing this cell-interface flux with the general expression  \eqref{ConvectiveDiffusiveFlux}, the numerical diffusive flux can be obtained as 
\begin{align}
d_I =  \frac{\left| s_I \right|}{2} \Delta U
\end{align}

It can be observed from \eqref{RHCondition} that $\Delta F$ and $\Delta U$ are $n \times 1$ vectors and thus a suitable choice for $s$ is an $n \times n$ matrix. One of the simplest assumptions for the matrix $s$ that can be conceived is a diagonal matrix with $n$ diagonal elements. Using this strategy the relation for obtaining the coefficient of numerical diffusion can be written as in \eqref{diagRHCondition}.  

\begin{equation}
\label{diagRHCondition}
 \Delta F_i = s_i \Delta U_i, \quad i= 1,2,\cdots n
\end{equation}

Thus the coefficient of numerical diffusion can be obtained as  \eqref{alphai}
\begin{align}
\label{alphai}
\alpha_{I,i} = \left| s_i \right| =\left| \frac{\Delta F_i}{\Delta U_i} \right|, \quad i =1,2,3 
\end{align} 
As $\Delta F_{i} = 0$ for stationary discontinuities, the numerical diffusion then vanishes, leading to exact capturing of grid-aligned discontinuities. It can be observed from (\ref{alphai}) that the coefficient $\alpha_I$ can go out of bounds when the denominator becomes small. 
\begin{equation}
 \Delta U \rightarrow 0,  ~~\alpha_I \rightarrow \infty.
\end{equation}
In order to introduce boundedness and stabilize the numerical scheme, 
$\alpha_I$ is to be restricted to a physically feasible range of eigenvalues of the flux Jacobian matrix. This process, termed as wave speed correction,  \eqref{wavespeedcorrection} is incorporated such that the coefficient of numerical diffusion lies within the eigenspecturm of the flux Jacobian \textit{i.e.}, $\alpha_I \in 
\left[\lambda_{max}, \lambda_{min}\right]$.
\begin{align}
\label{wavespeedcorrection}
 \lvert \alpha_{\mathrm{I}} \rvert = \begin{cases}
                             \lambda_{\textit{max}},~~ if ~~ \lvert \alpha_{\mathrm{I}} \rvert > \lambda_{\textit{max}} \\
                             \lambda_{\textit{min}},~~ if ~~ \lvert \alpha_{\mathrm{I}} \rvert < \lambda_{\textit{min}}\\
                             \lvert \alpha_{\mathrm{I}} \rvert , ~~ \textrm{otherwise} \\
                            \end{cases}\\
                            \label{FinalInterfaceflux}
 F_{\mathrm{I}} = \frac{1}{2}\left[ F_L + F_R\right] - \frac{\lvert\alpha_{\mathrm{I}}\rvert}{2}\left[ U_{R} - U_{L} \right].
\end{align}
Hence the final numerical flux at the cell-interface in MOVERS is given by \eqref{alphai}, \eqref{wavespeedcorrection} and \eqref{FinalInterfaceflux}.  This method is independent of eigen-structure of the underlying hyperbolic systems, is simple and can capture grid-aligned stationary discontinuities exactly. Two variations of MOVERS are introduced in  \cite{Jaisankar_SVRRao}:  $(i)$ an $n$-wave based coefficient of numerical diffusion, corresponding to $n$ number of conservation laws (MOVERS-n) and $(ii)$ a scalar diffusion, corresponding to the energy equation (as it contains the maximum of information), referred to as MOVERS-1. 
The robustness of the basic scheme has been improvised through its variants by Maruthi N.H. \cite{Maruthi_Thesis} and extended to other hyperbolic systems for magnetohydrodynamics and shallow water flows.  The simplicity and accuracy of MOVERS  makes this scheme a well-suited base-line solver for further research, apart from its independency of the eigenstrucure.  In this work this algorithm is chosen as the foundation to devise a new and efficient algorithm, named as MOVERS+. First, the wave-speed correction mechanism is removed by a reformulation of the basic Rankine-Hugoniot solver.  Further, exact schock capturing is deliberately given up for enhancing robustness but exact contact discontinuity capturing is retained for preserving accuracy.    

\subsection{A new central solver: \emph{MOVERS+}} \label{sec:movers+}
MOVERS \cite{Jaisankar_SVRRao}   requires wave speed correction in order to restrict the coefficient of diffusion to within the eigenspectrum. To avoid wave speed correction, a simpler strategy is proposed in this section which is described below.
\begin{align}
\label{alphai1}
d_{\mathrm{I,j}} &=  \frac{1}{2}\left| \frac{\Delta F_j}{\Delta U_j}\right| \Delta U_j, \quad j =1,2,3 \\
	 &= \frac{1}{2}\frac{\left |\Delta F_j \right |}{\sign(\Delta U_j) \Delta U_j } \Delta U_j  \\ 
	 & = \frac{1}{2}\sign(\Delta U_j)\left| \Delta F_j \right| , \quad j =1,2,3
\end{align}
where the relation $\frac{1}{\sign(\cdot)} = \sign(\cdot)$ is used.  
This form of $d_I$ will eliminate the need of wave speed correction for MOVERS. Numerical experimentation has revealed that this numerical scheme has very low diffusion and captures steady discontinuities exactly but encounters problems in smooth regions due to lack of sufficient numerical diffusion. Therefore, using a shock sensor (\ref{ShockSensor}), an additional numerical diffusion is introduced. This additional diffusion is based on the fluid velocity, which is demonstrated to be sufficient to avoid unphysical expansions in smooth regions \cite{Venkat_thesis,Ramesh_Thesis}.  
The coefficient of numerical diffusion for MOVERS+ is given by
\begin{align}
\label{MOVERS_NWSC}
 \lvert d_{\mathrm{I}} \rvert_j =  \Phi Sign(\Delta U_j)\lvert \Delta F_j \rvert + \left(\frac{|V_{nL}|+|V_{nR}|}{2}\right) \Delta U_j, \quad j =1,2,3 
\end{align} 
where the $\Phi$ is the shock sensor defined by  
\begin{align}
\label{ShockSensor}
\Phi = \left| \frac{\Delta p}{2 {p_{\mathrm{I}}}}\right| \ \textrm{with} \ p_{I} = \frac{p_L+p_R}{2} 
\end{align}

The features of this modified algorithm, MOVERS+, are as follows.  
\begin{itemize}
 \item It can capture steady grid-aligned contact discontinuities exactly and provides low diffusion otherwise.  
 \item It has sufficient numerical diffusion near shocks so as to avoid shock instabilities (deliberately giving up exact shock capturing of MOVERS for gain in robustness).  
 \item It does not need entropy fix for smooth regions or in expansion regions.
 \item It does not require any wave speed correction, unlike in MOVERS. 
 \item It is a simple central solver and is not based on Riemann solvers, field-by-field decompositions or  complicated flux splittings, thus making it a suitable candidate for further extensions.  
\end{itemize}  

\section{Results and Discussion}\label{sec:results}
To test the accuracy and robustness of the numerical schemes RICCA and MOVERS+, results from various 1-D and 2-D benchmark cases are presented in the following sections.
\subsection{1-D Shock tube problems}
These are robustness and accuracy test problems with initial conditions as specified in chapter 4 of \cite{TORO_1} and \cite{Quirk}. Both the schemes RICCA and MOVERS+ are tested for: Sod test problem with a sonic point, a strong shock, strong rarefaction, slowly-moving shock and slowly-moving contact-discontinuity and their interactions.
The initial conditions for these test cases are given in the Table \ref{1D_Testcases_Table}. For all the test cases that are being considered in 1D a total of 100 computation cells are considered and the CFL number, unless and otherwise specified, is taken as 0.1. Numerical results are compared with the analytical solutions of the Riemann problems.  
\begin{table}[h!]
\centering\begin{large}          \end{large}
\begin{tabular}{ |c|c|c|c|c|c|c|c|c|}
\hline
Case & $\rho_L$&    $p_L$     &     $u_L$     &     $\rho_R    $     &     $p_R$         &     $u_R$  \\
\hline
1         &1.0    &    1.0    &     0.0         &    0.125        &    0.1        &     0.0  \\
\hline
2        &1.0    &    0.4    &    -2.0        &    1.0            &    0.4        &    2.0 \\
\hline
3        &1.0    &    1000.0&    0.0        &    1.0            &    0.01        &    0.0 \\
\hline
4        &1.0    &    0.01    &    0.0        &    1.0            &    100.0    &    0.0\\
\hline
5        &5.99924    &460.894    &19.5975    &    5.99242        &    46.0950    &    -6.19633\\
\hline
6        &1.0    &    $\frac{1}{\gamma M^2}$    &    1.0    & $\frac{\frac{\gamma + 1}{\gamma -1}\frac{p_R}{p_L}+1}{\frac{\gamma -1}{\gamma + 1}+\frac{p_R}{p_L}}$&$P_L\frac{2\gamma M^2 - (\gamma - 1)}{\gamma +1}$&    $\sqrt(\gamma \frac{(2 + (\gamma -1)M^2)p_R}{(2\gamma M^2 + (1 - \gamma)) \rho_R})$\\
\hline
7        &1.4    &    0.4    &    0.0        &    1.0            &    0.4        &    0.0 \\
\hline
8        &1.4    &    1.0    &    0.1        &    1.0            &    1.0        &    0.1\\
\hline
9        &3.86&    10.33&    -0.81        &    1.0            &    1.0        &    -3.44\\
\hline
\end{tabular}
\caption{Initial conditions for 1D test cases as given in \cite{TORO_1}}
\label{1D_Testcases_Table}
\end{table} \\

\begin{figure}[h!]
\begin{subfigure}[t]{0.475\textwidth}
\begin{center}
\includegraphics[width=\textwidth]{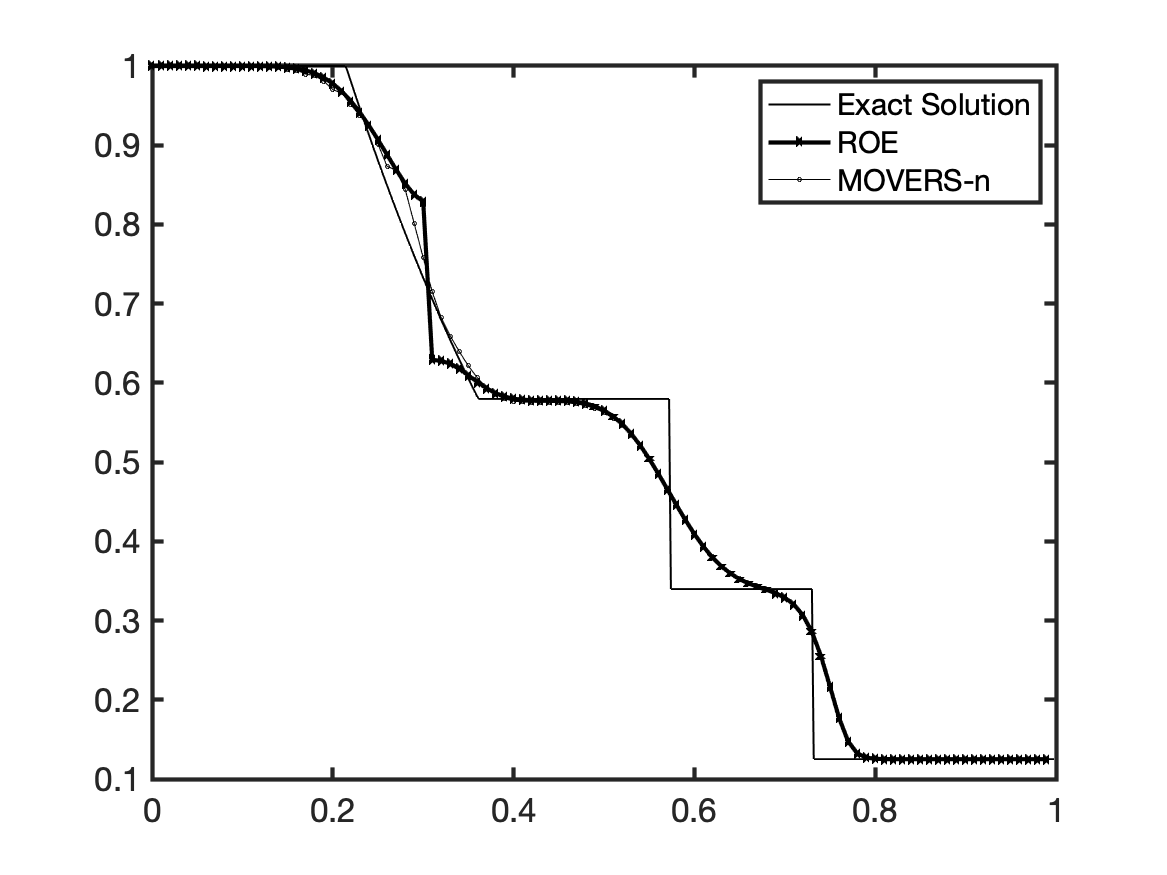}
\caption{Non-smoothness in expansion and expansion shock}
\label{Sonic_Glitch}
\end{center}
\end{subfigure}
\begin{subfigure}[t]{0.475\textwidth}
\begin{center}
\includegraphics[width=\textwidth]{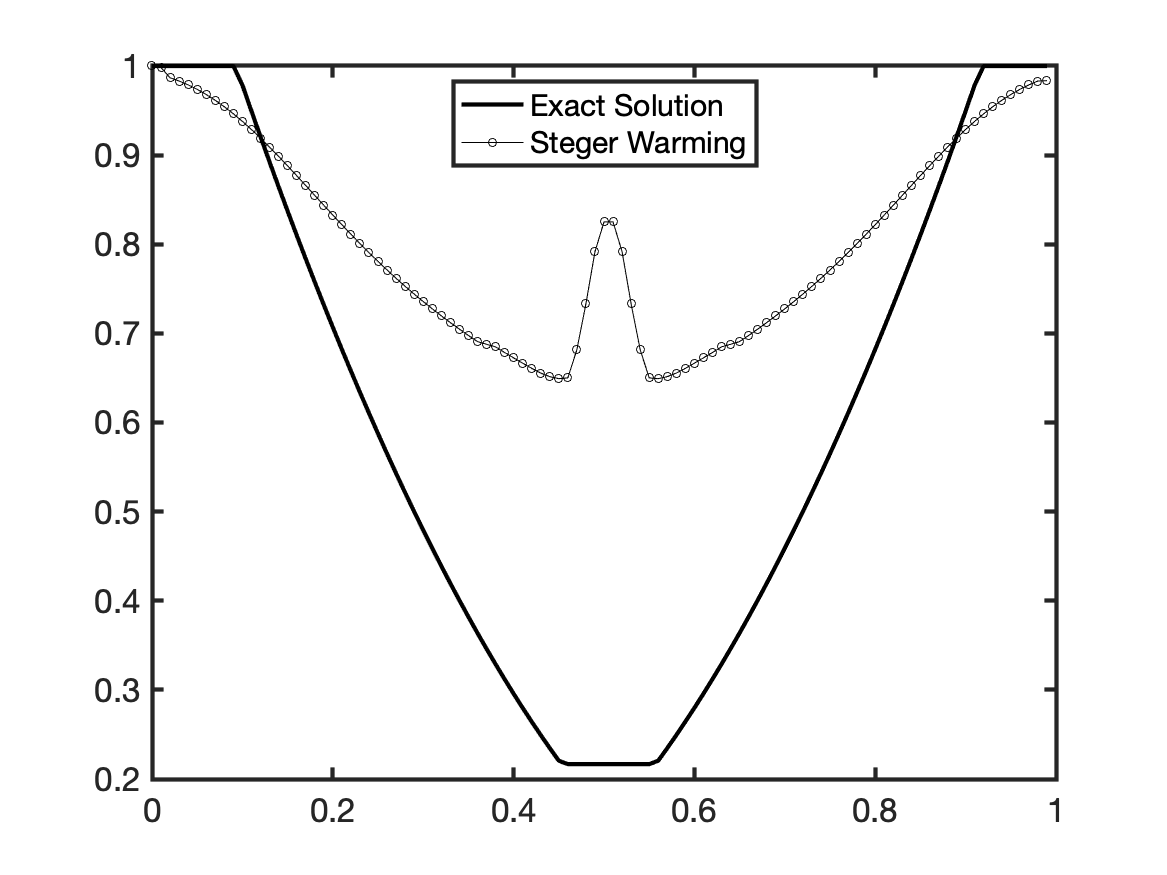}
\caption{Internal energy plot for Toro test case 2 }
\label{Test2_Internalenergy}
\end{center}
\end{subfigure}
\begin{subfigure}[t]{0.475\textwidth}
\includegraphics[width=\textwidth]{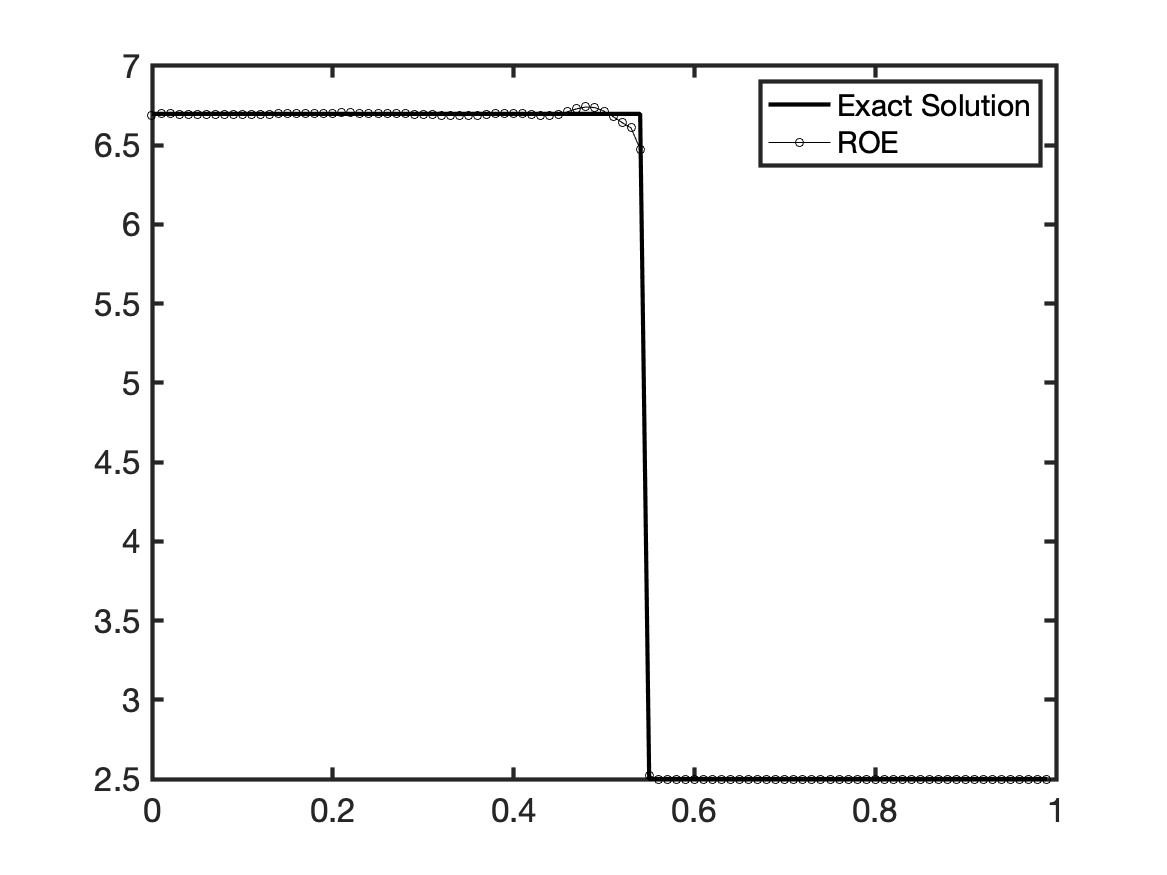}
\caption{Post shock oscillations in Slowly moving shock}
\label{Slowly_Moving_Shock_Discontinuity_ROE}
\end{subfigure}
\caption{General Issues with popular schemes}
\end{figure}

\begin{figure}[h!]
\begin{center}
\includegraphics[scale=0.8]{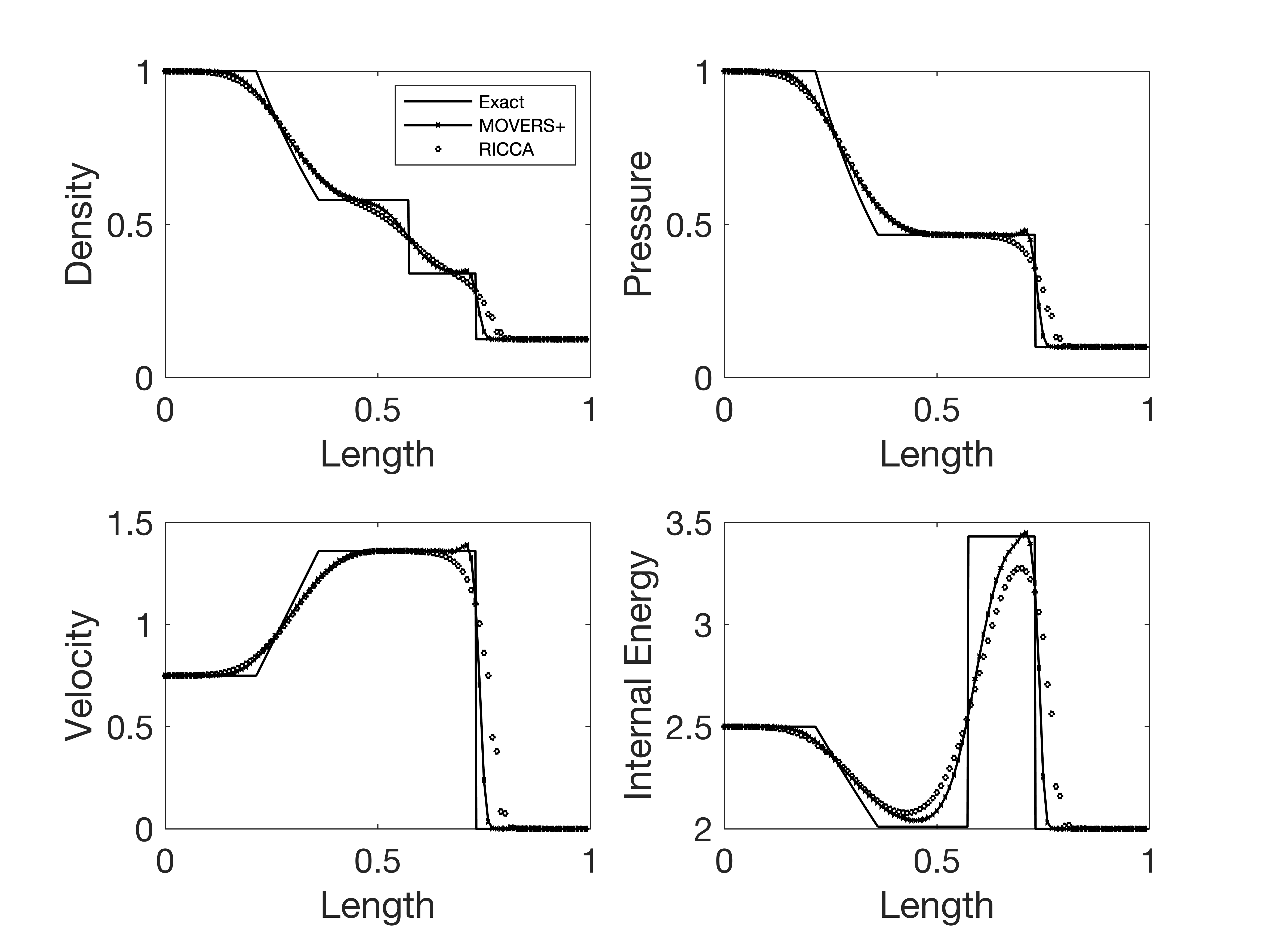}
\caption{Results for Toro test case 1}
\label{Plus_Toro_Test_Case1}
\end{center}
\end{figure}

Test case 1 corresponds to a Sod shock tube problem with mild shock strength. This test case has an expansion fan (containing a sonic point) moving to the left, a shock moving towards right side and a contact-discontinuity in between these two. Typically low diffusion schemes encounter problems in the expansion fans, especially at sonic points. MOVERS-n gives a small non-smooth variation near the sonic point while Roe scheme yields a large and unphysical expansion shock, as shown in figure (\ref{Sonic_Glitch}).
It can be observed from figure  (\ref{Plus_Toro_Test_Case1}) that both RICCA and MOVERS+ do not produce expansion shocks or non-smoothness in the expansion region. Further, It can be observed that RICCA is more diffusive in shock capturing when compared to MOVERS+, which is also seen in many other test cases described in this section.

\begin{figure}[h!]
\begin{center}
\includegraphics[scale=0.8]{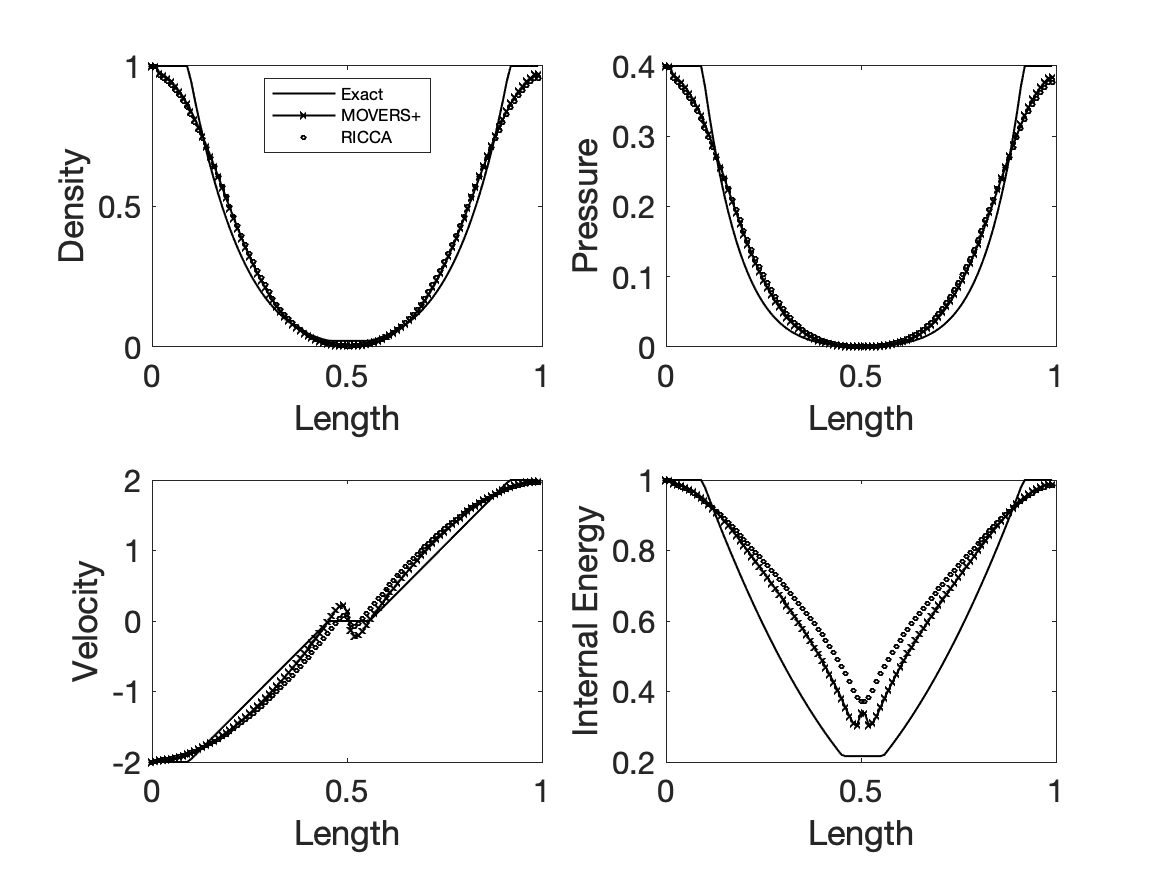}
\caption{Results for Toro test case 2}
\label{Plus_Toro_Test_Case2}
\end{center}
\end{figure}

\begin{figure}[h!]
\begin{center}
\includegraphics[scale=0.8]{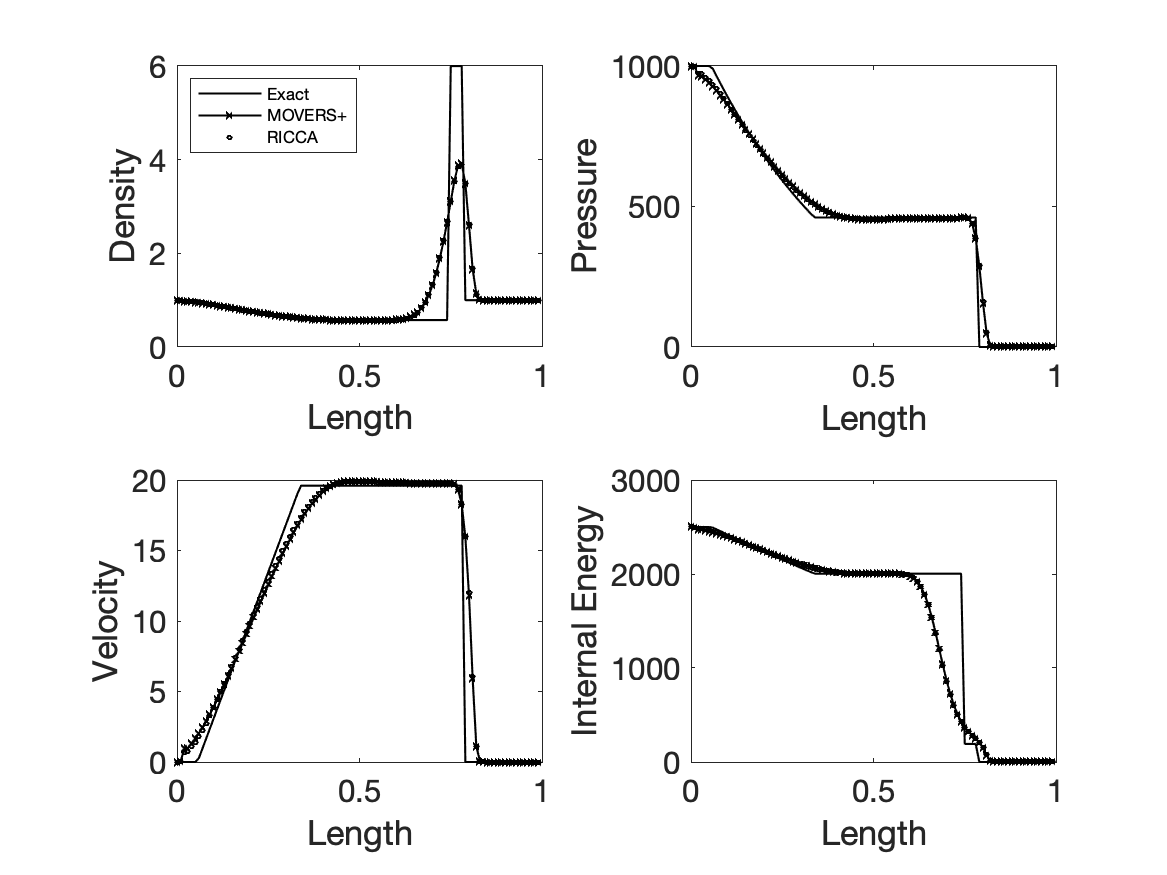}
\caption{Results for Toro test case 3}
\label{Plus_Toro_Test_Case3}
\end{center}
\end{figure}

Test case 2, also known as \emph{123 problem}, consists of two strong and symmetric rarefactions approaching each other and a trivial contact-discontinuity of zero wave speed between them.  The pressure between these rarefactions can go as low as zero (close to vacuum).  This test case can also be considered to be a benchmark test case for low density flows. Many low diffusion schemes fail for this case. Other numerical schemes would give a wrong value of internal energy as shown in figure (\ref{Test2_Internalenergy}). Hence the schemes which can resolve this test case are considered to be robust. 
Figure (\ref{Plus_Toro_Test_Case2}) shows the results obtained by MOVERS+ and RICCA respectively. The ability of a numerical scheme to capture the low density is tested in this test case.  Both the new numerical schemes can capture the low pressure and density regions close to vacuum.
Test case 3 represents a strong and severe problem, designed to test the robustness and accuracy of numerical schemes. It consists of a left rarefaction, a contact-discontinuity and a strong right shock wave with shock Mach number 198. This test case forms the left half portion of the blast wave problem of Woodward and Colella \cite{Woodward}.  Figure (\ref{Plus_Toro_Test_Case3}) represents  solutions obtained from MOVERS+ and RICCA. Both the numerical schemes are capable of capturing strong shocks.

\begin{figure}[h!]
\begin{center}
\includegraphics[scale=0.8]{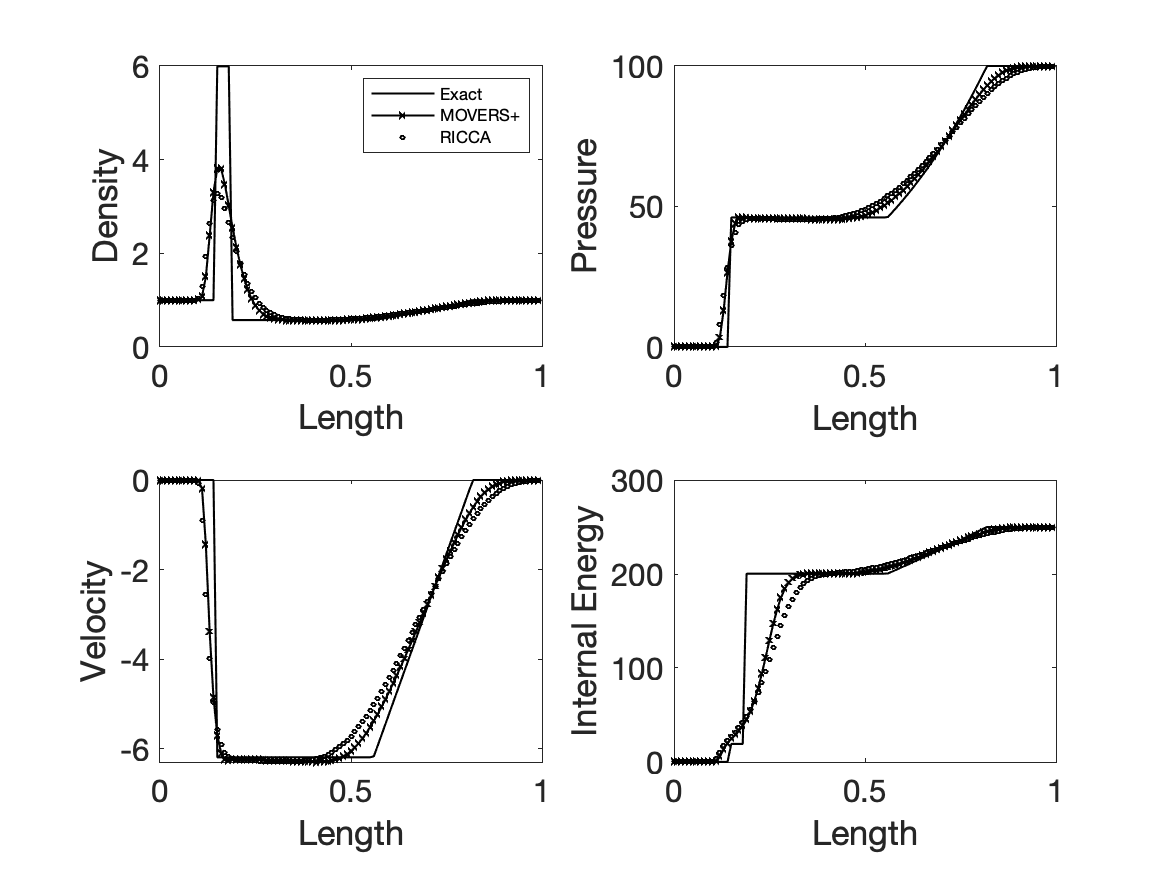}
\caption{Results for Toro test case 4}
\label{Plus_Toro_Test_Case4}
\end{center}
\end{figure}

Test case 4 represents another difficult problem, which is taken from the right half of the Woodward and Colella \cite{Woodward} problem. This problem consists of a left shock, a contact-discontinuity and a right rarefaction. Figures (\ref{Plus_Toro_Test_Case4}) represent the solutions of the test case 4 using MOVERS+ and RICCA respectively. 
%
%
%
%

\begin{figure}[h!]
\begin{center}
\includegraphics[scale=0.8]{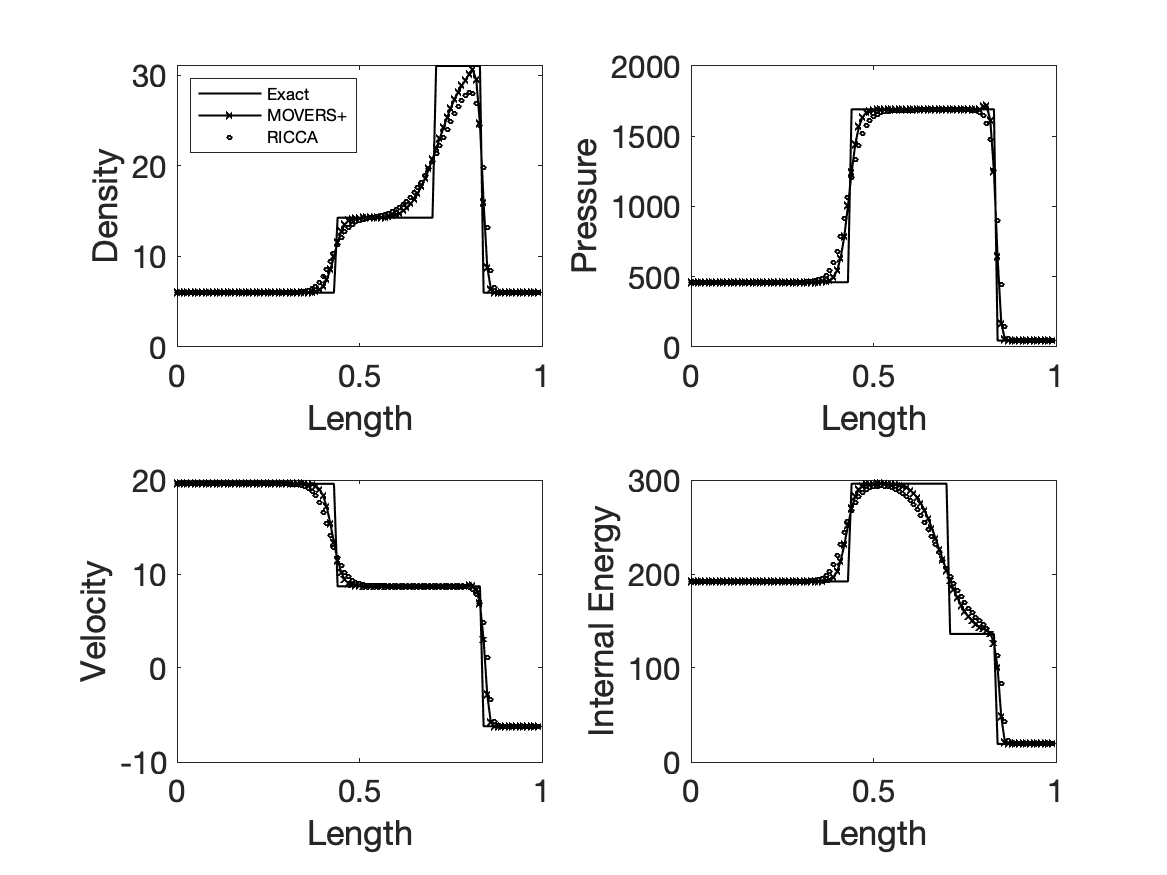}
\caption{Results for Toro test case 5}
\label{Plus_Toro_Test_Case5}
\end{center}
\end{figure}

Test case 5 is made up of solutions emerging from test cases 3 and 4 and represents the wave system resulting from the interaction of two strong shock waves propagating in opposite directions.  The solution consists of a slowly moving left shock, a contact-discontinuity and a right travelling shock wave. Figure (\ref{Plus_Toro_Test_Case5}) represents the solutions obtained using MOVERS+ and RICCA.

\begin{figure}[h!]
\begin{center}
\includegraphics[scale=0.8]{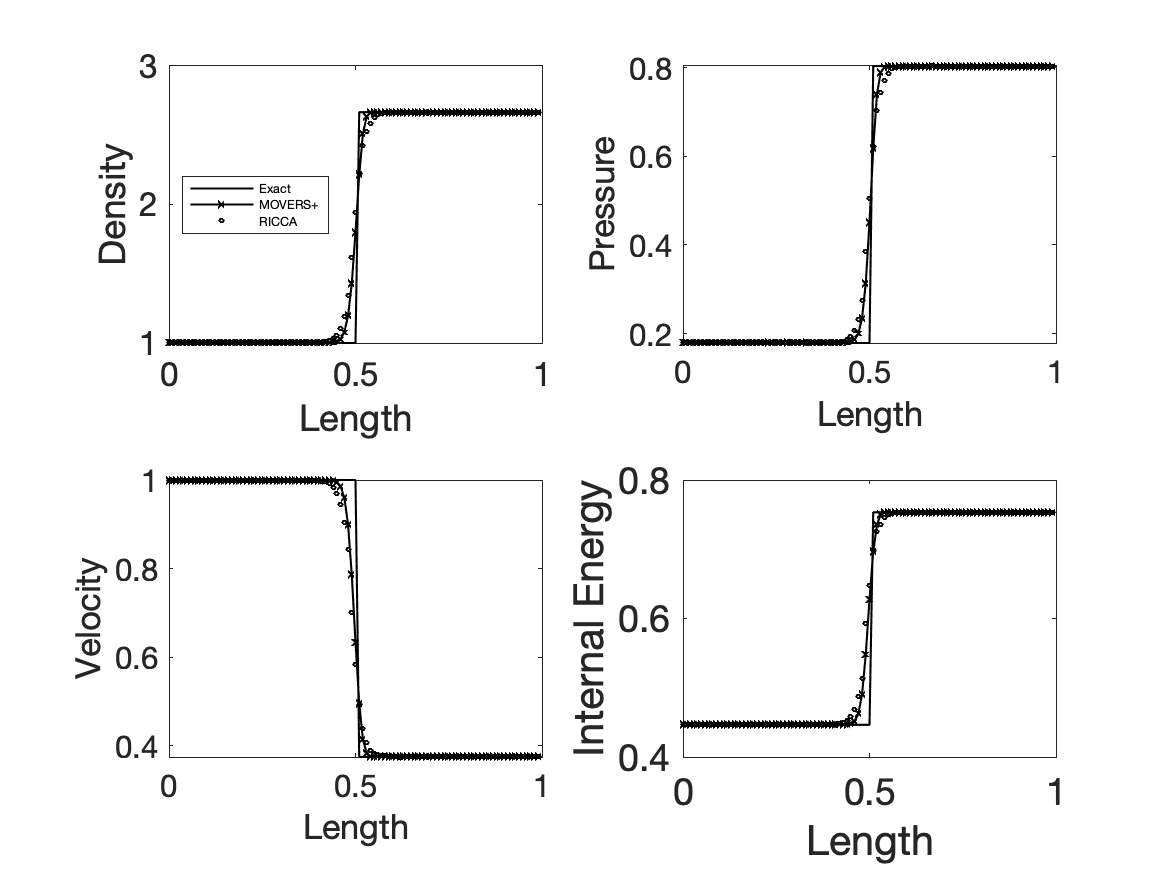}
\caption{Results for steady shock test case}
\label{Steady_Shock_Plus}
\end{center}
\end{figure}

Test case 6 is designed to mimic the conditions of shock in steady state \cite{Zhang_Shu}. Initial conditions for the this test case are given in the table (\ref{1D_Testcases_Table}). MOVERS-n and Roe schemes capture steady shock exactly. Figure (\ref{Steady_Shock_Plus}) represents the results for steady state shock using MOVERS+ and RICCA.   Both the numerical schemes diffuse the steady shock. This is expected, as the exact shock capturing is deliberately given up in the designing of the schemes for avoiding shock instabilities, still retaining exact contact discontinuity capturing.   

\begin{figure}[h!]
\begin{center}
\includegraphics[scale=0.8]{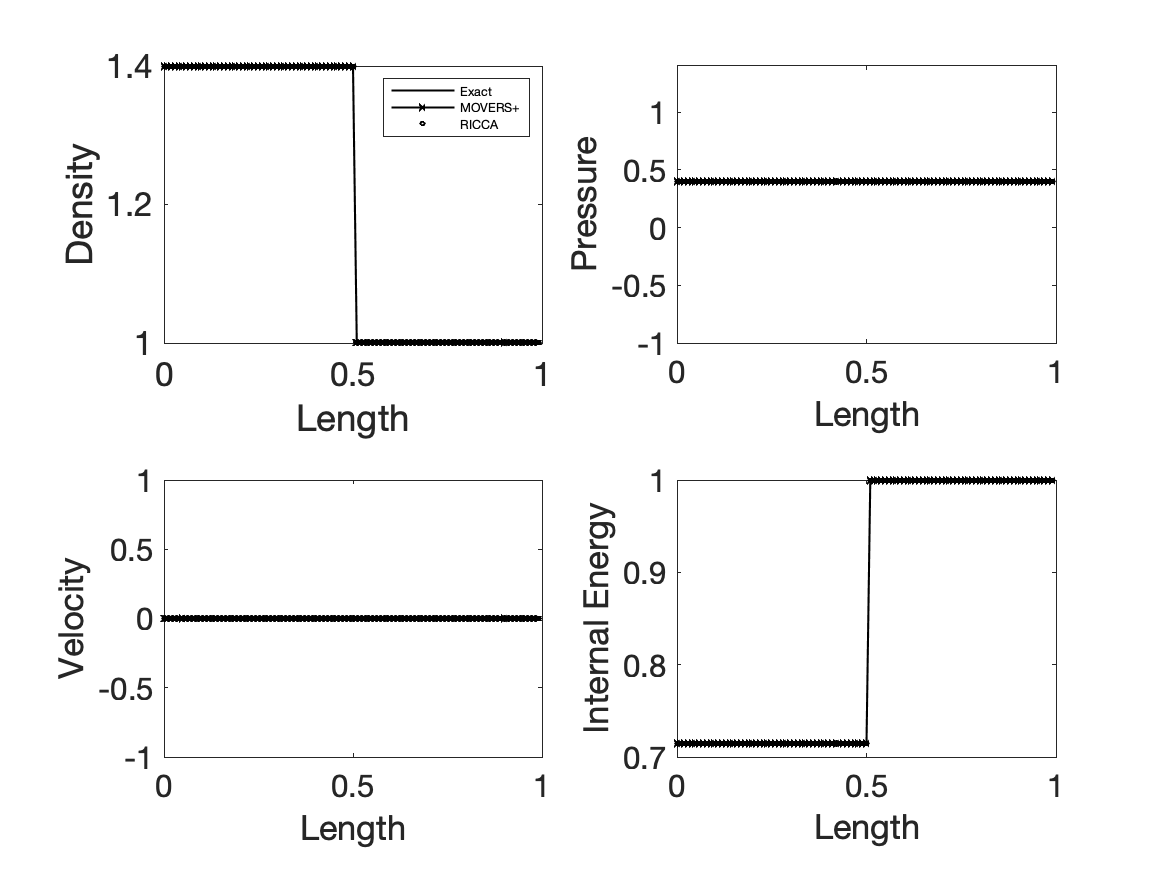}
\caption{Results for steady contact-discontinuity test case}
\label{Steady_Contact}
\end{center}
\end{figure}

\begin{figure}[h!]
\begin{center}
\includegraphics[scale=0.8]{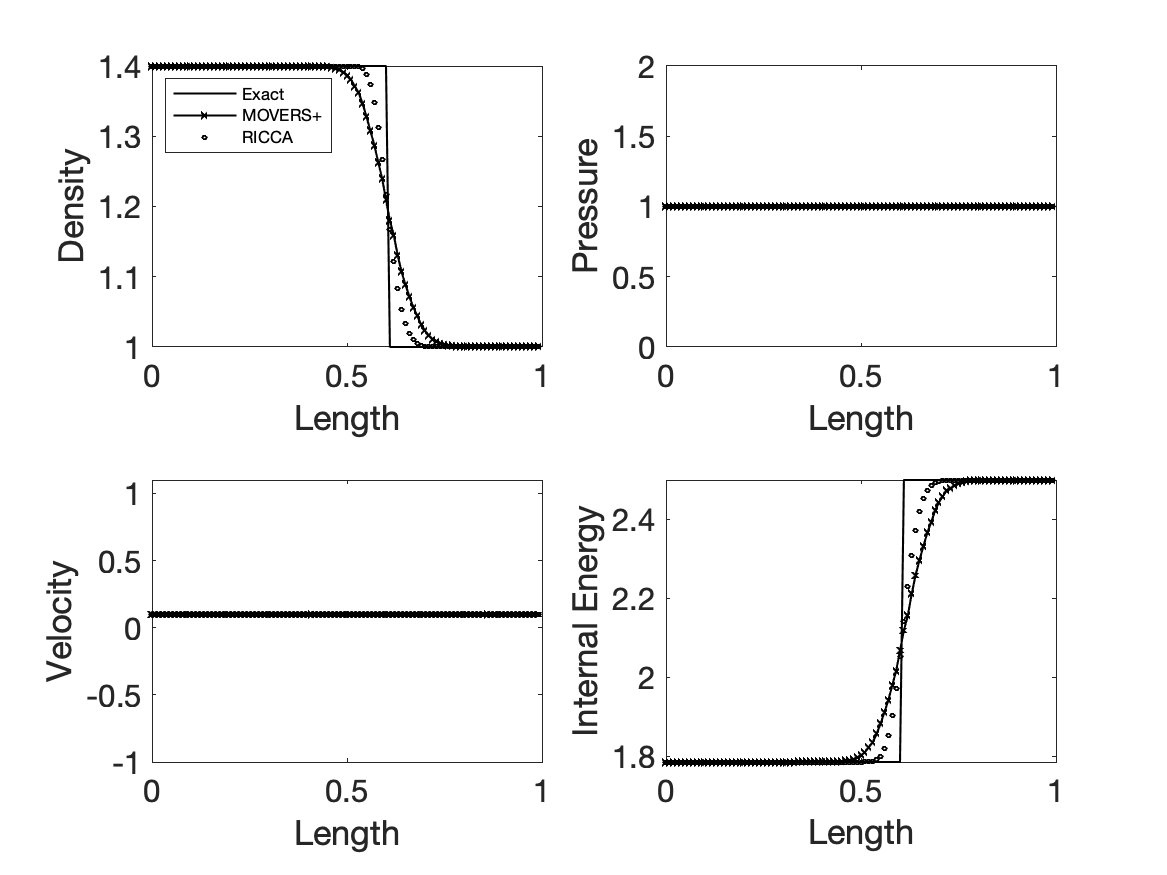}
\caption{Results for slowly moving contact-discontinuity test case}
\label{Slowly_Moving_Contact_Discontinuity_Plus}
\end{center}
\end{figure}

\begin{figure}[h!]
\begin{center}
\includegraphics[scale=0.8]{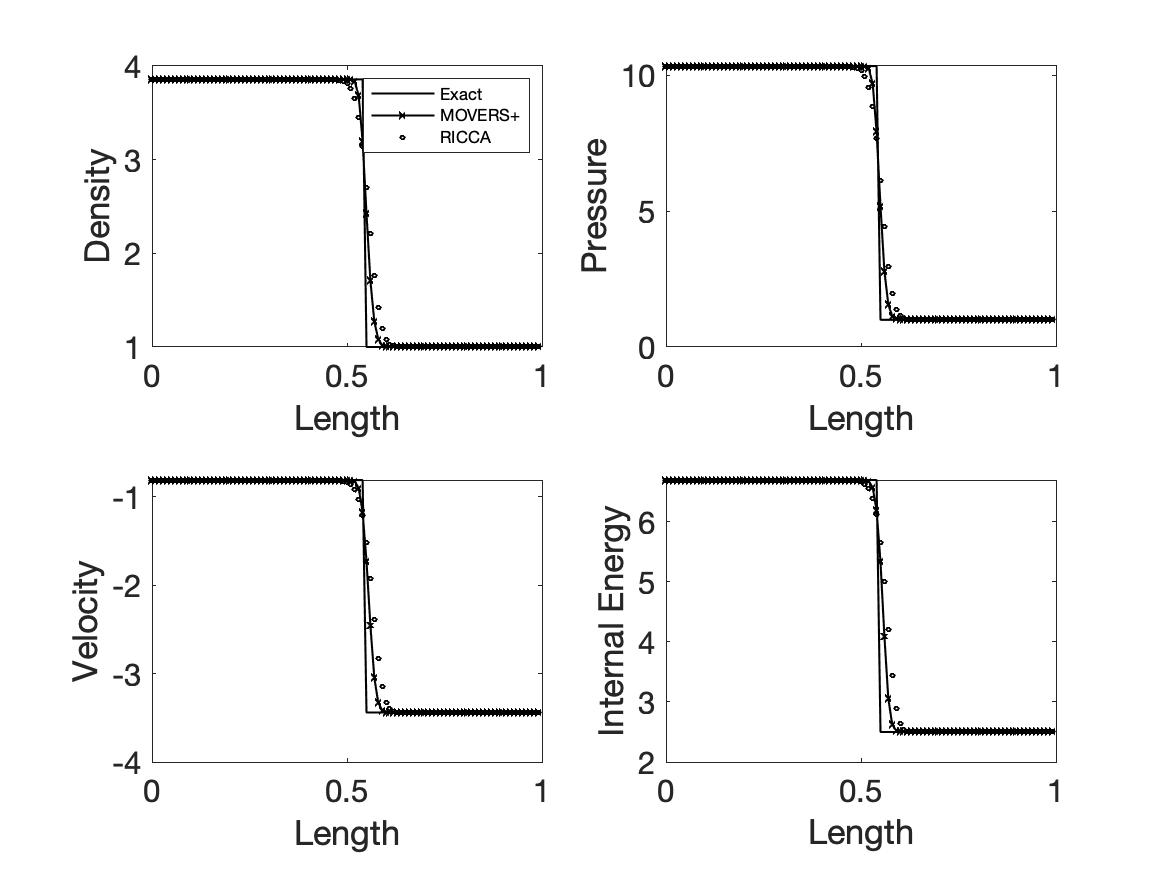}
\caption{ Results for slowly moving shock test case}
\label{RICCA_Slowly_Moving_Shock_Discontinuties}
\end{center}
\end{figure}

 Both numerical schemes are designed to capture exactly the steady state contact-discontinuities. In order to check this capability, test case 7 is designed to specifically mimic a steady state contact-discontinuity. It is a known fact that across this discontinuity there will be no variation in pressure and velocity but density variation occurs. Both these schemes have the ability to capture steady contact-discontinuity exactly as shown in figure (\ref{Steady_Contact}).
 Typically low diffusive schemes will generate oscillations near slowly moving shocks and contact discontinuities, as shown in (\ref{Slowly_Moving_Shock_Discontinuity_ROE}) and as described in  \cite{Quirk,shijin,stiriba,karni2}.  Figures (\ref{Slowly_Moving_Contact_Discontinuity_Plus}) and (\ref{RICCA_Slowly_Moving_Shock_Discontinuties}) represent solutions obtained by RICCA and MOVERS+ which show no such oscillations.  
\subsection{2-D Euler test cases}
In this section a set of 2D benchmark test cases for Euler equations are considered to check the accuracy and robustness of RICCA and MOVERS+. 

\subsubsection{Oblique shock reflection} 

\begin{figure}[h!]
\makebox[\textwidth][c]{%
\begin{subfigure}[b]{.5\textwidth}
\centering
\includegraphics[width=\textwidth,height=\textheight,keepaspectratio]{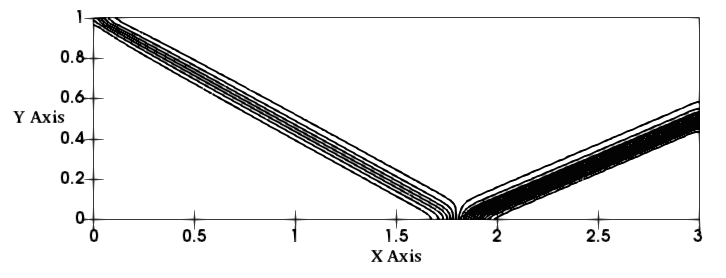}
\caption{Second order accurate on $120\times40$ grid}
\end{subfigure}%
\begin{subfigure}[b]{.5\textwidth}
\centering
\includegraphics[width=\textwidth,height=\textheight,keepaspectratio]{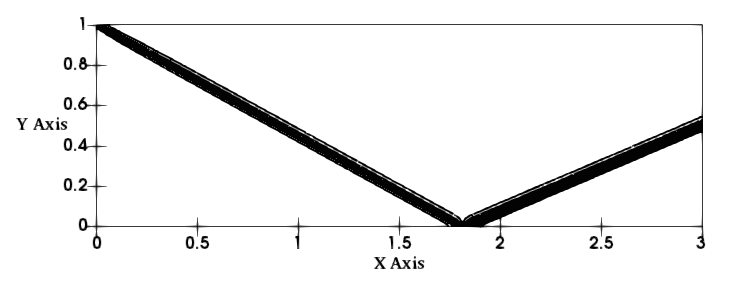}
\caption{Second order accurate on $240\times80$}
\end{subfigure}
}
\caption{MOVERS+ - Pressure Contours (0.71:0.1:2.91) - for regular shock reflection}
\label{fig:NWSC_Regular_shock_reflection_FO_SO_120_40_ch3}
\end{figure}

\begin{figure}[h!]
\makebox[\textwidth][c]{%
\begin{subfigure}[b]{.5\textwidth}
\centering
\includegraphics[width=\textwidth,height=\textheight,keepaspectratio]{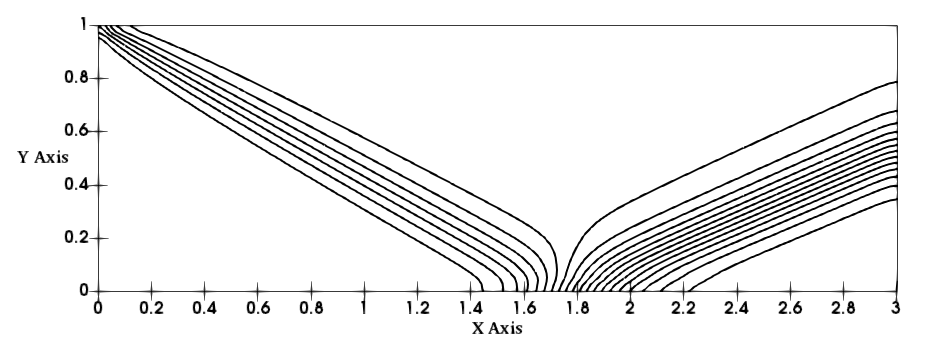}
\caption{Second order accurate on $120\times40$ grid}
\end{subfigure}%
\begin{subfigure}[b]{.5\textwidth}
\centering
\includegraphics[width=\textwidth,height=\textheight,keepaspectratio]{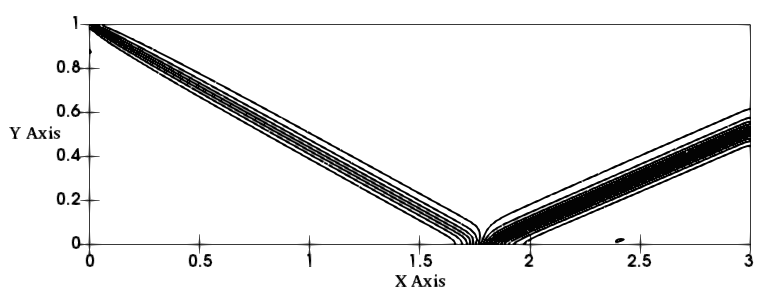}
\caption{Second order accurate on $240\times80$ grid}
\end{subfigure}
}
\caption{RICCA - Pressure Contours (0.71:0.1:2.91) - for regular shock reflection}
\label{fig:RICCA_Regular_shock_reflection_FO_SO_120_40_ch3}
\end{figure}

\begin{figure}[h!]
\makebox[\textwidth][c]{%
\begin{subfigure}[b]{.5\textwidth}
\centering
\includegraphics[width=\textwidth,height=\textheight,keepaspectratio]{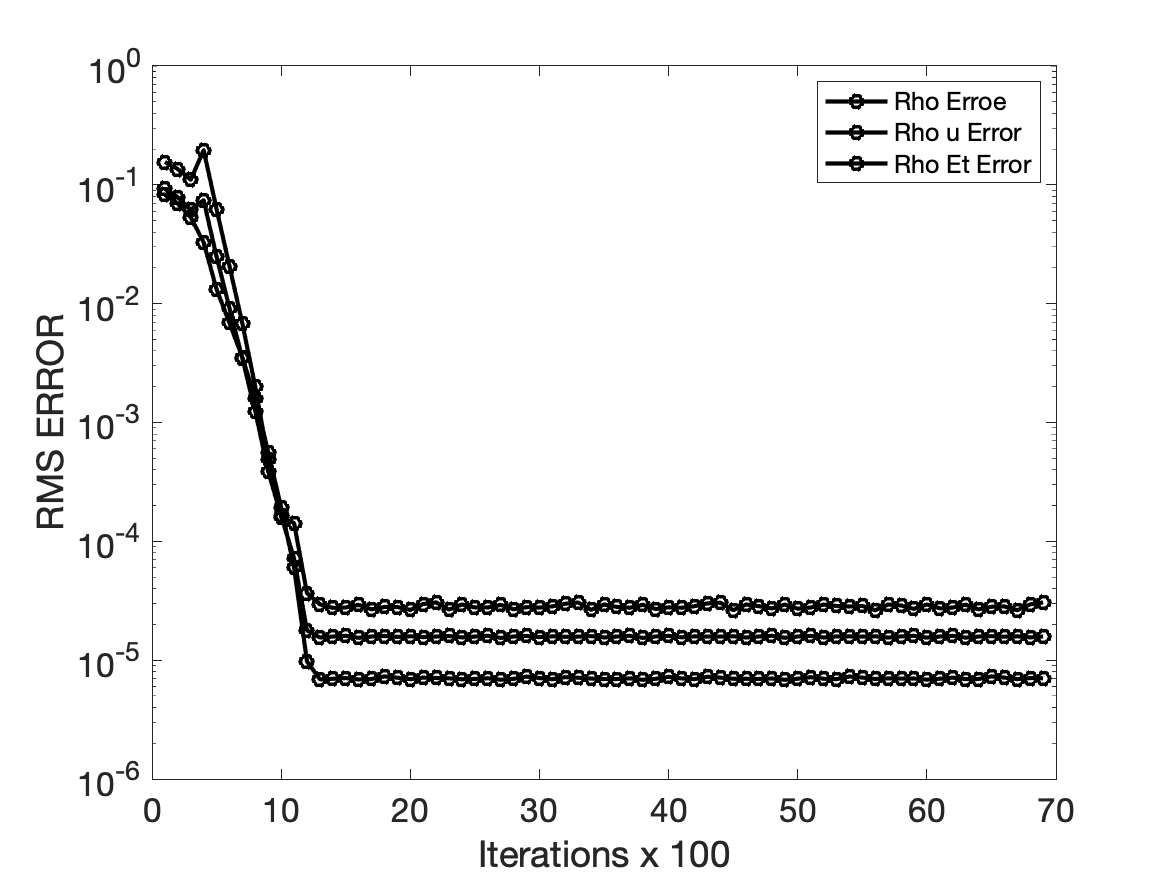}
\caption{MOVERS+}
\end{subfigure}
\begin{subfigure}[b]{.5\textwidth}
\centering
\includegraphics[width=\textwidth,height=\textheight,keepaspectratio]{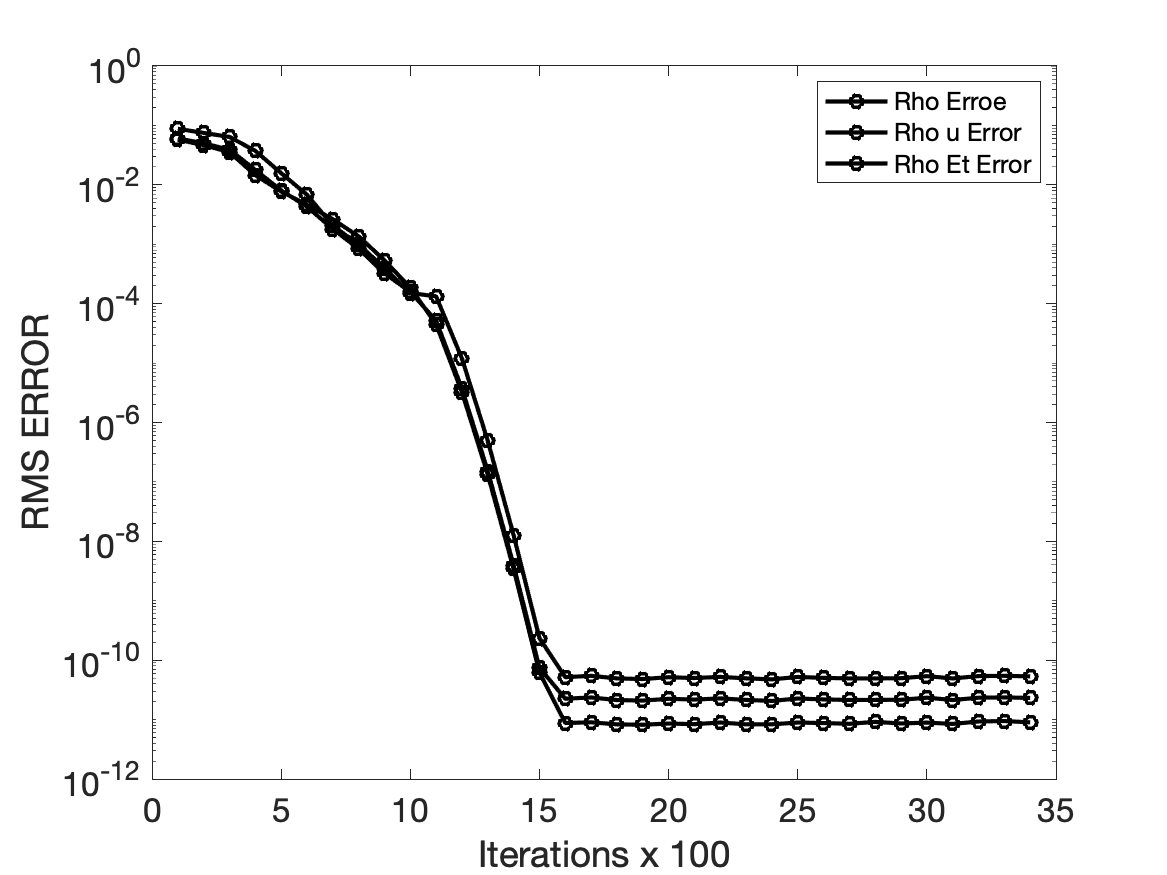}
\caption{RICCA}
\end{subfigure}%
}
\caption{Convergence plot of MOVERS+ and RICCA on $240 \times 80$ grid}
\label{convergnceplotreflectionproblem_RICCA}
\end{figure}

This test case~\cite{Yee} describes an oblique shock hitting a flat plate and getting reflected. The computational domain considered for this test case is $\left[0, 3\right] \times \left[0, 1\right]$. An oblique shock with the incident shock angle of $29^{\circ}$ and the free stream Mach number $M = 2.9$ is introduced from the left-top corner of the computational domain.   The initial conditions for this test problem are as follows.  
\begin{align*}
(\rho,u,v,p)_{0,y,t} &=(1.0,2.9,0,1/1.4), ~~ \textrm{(Inflow Conditions)}\\
(\rho,u,v,p)_{x,1,t} &=(1.69997,2.61934,-0.50633,1.52819),~~ \textrm{(Post Shock Conditions)}
\end{align*}
with $(\cdot)_{0,y,t}$ corresponding to the conditions given on the left side boundary and $(\cdot)_{x,1,t}$ on the top boundary.  
Flow tangency boundary conditions are imposed at the wall boundary, which is the bottom part of the computational domain, and supersonic outflow boundary conditions are used at the right side of the computational domain. Standard grid size of $120 \times 40 $ and $240 \times 80$ are considered for the study. Since the flow features considered are steady-state conditions the code is run till the relative error reaches machine epsilon or the number of iterations reaches 100000. Figures (\ref{fig:NWSC_Regular_shock_reflection_FO_SO_120_40_ch3}) and  (\ref{fig:RICCA_Regular_shock_reflection_FO_SO_120_40_ch3}) show results with second-order accurate versions of the schemes for the oblique shock reflection problem.  Both the schemes,  RICCA and MOVERS+, are capable of resolving incident and reflected shocks well.  Further, it can be observed that the shock resolution capability of MOVERS+ is much better when compared to RICCA on a coarse grid.  The convergence plots of second order MOVERS+ and RICCA are shown in figure (\ref{convergnceplotreflectionproblem_RICCA}).  It can be observed that RICCA converges to machine epsilon whereas MOVERS+ converges to $1\times 10^{-7}$.

\subsubsection{Supersonic flow over $15^o$ compression ramp}
\begin{figure}[h!]
\makebox[\textwidth][c]{%
\begin{subfigure}[b]{.5\textwidth}
\centering
\includegraphics[width=\linewidth]{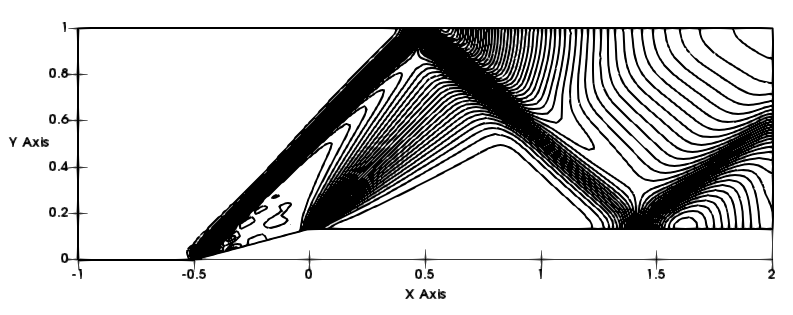}
\caption{MOVERS+ on $121\times 41$ Grid}
\end{subfigure}
\begin{subfigure}[b]{.5\textwidth}
\centering
\includegraphics[width=\linewidth]{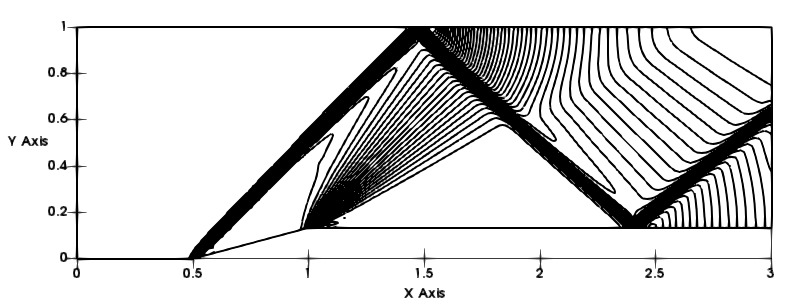}
\caption{MOVERS+ on $241\times 81$ Grid}
\end{subfigure}%
}
\caption{Second order results of Mach 2 flow on $15^o$ ramp using MOVERS+}
\label{fig:MOVERS+_RampFlow}
\end{figure}

\begin{figure}[h!]
\makebox[\textwidth][c]{%
\begin{subfigure}[b]{.5\textwidth}
\centering
\includegraphics[width=\linewidth]{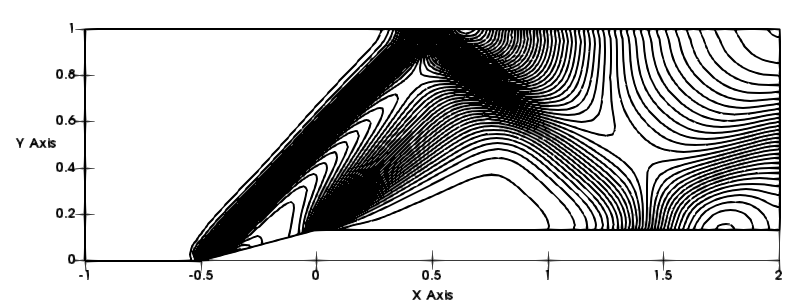}
\caption{RICCA on $121\times 41$ Grid}
\end{subfigure}
\begin{subfigure}[b]{.5\textwidth}
\centering
\includegraphics[width=\linewidth]{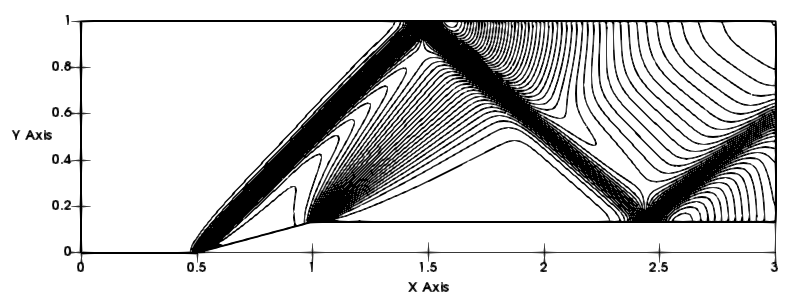}
\caption{RICCA on $241 \times 81$ Grid}
\end{subfigure}
}
\caption{Second order results of Mach 2 flow on $15^o$ ramp using RICCA}
\label{fig:RICCA_RampFlow}
\end{figure}

\begin{figure}[h!]
\makebox[\textwidth][c]{%
\begin{subfigure}[b]{.5\textwidth}
\centering
\includegraphics[width=\linewidth]{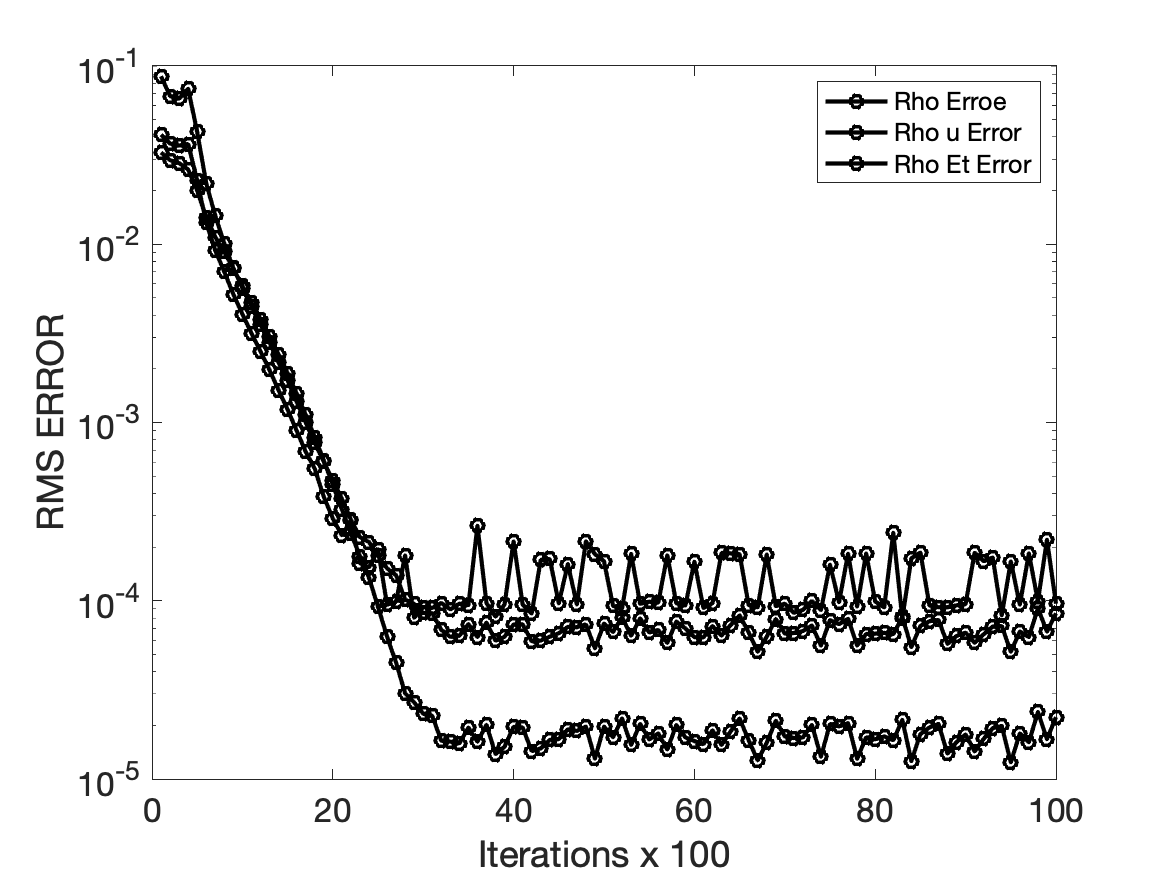}
\caption{MOVERS+ residual on $240\times 80$ Grid}
\end{subfigure}
\begin{subfigure}[b]{.5\textwidth}
\centering
\includegraphics[width=\linewidth]{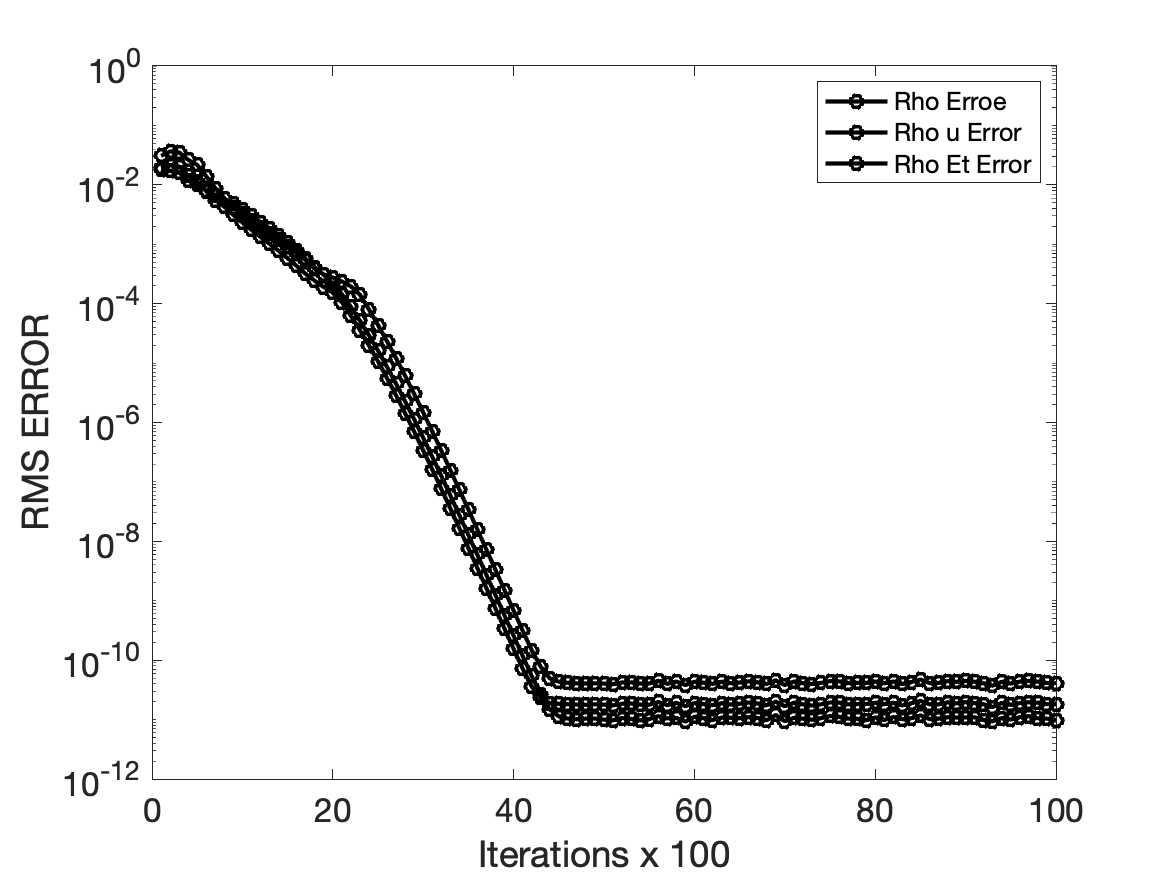}
\caption{RICCA  residual on $240\times 80$ Grid}
\end{subfigure}
}
\caption{Relative error for 2O simulations}
\label{fig:RICCA_MOVERSP_Convergence}
\end{figure}

This is the case of a supersonic flow over a $15^o$ degree ramp \cite{levy} placed in a wind tunnel. The computational domain consists of the following dimensions $[0,3]\times[0,1]$ with supersonic inlet, supersonic exit and flow tangency wall boundary conditions applied. As a supersonic flow with Mach number $M =2.0$ approaches this ramp, an oblique shock emerges from the lower end corner (compression corner) of the ramp and an expansion fan emerging from the upper corner (expansion corner) of the ramp.  The oblique shock emerging from the compression corner at the beginning of the ramp reflects from the upper wall.  The reflected shock interacts with the emerging expansion fan and further gets reflected from the bottom wall.  Following grid sizes are used in the numerical simulation: $120 \times 40$ and $240 \times 80$. Figures (\ref{fig:MOVERS+_RampFlow}) and (\ref{fig:RICCA_RampFlow}) show the second order accurate results from MOVERS+ and RICCA on two different grid sizes. It can be observed from these figures that the second order results of MOVERS+ on coarse grid are comparable with the second order results of RICCA on fine grid. The relative error for this test case is plotted for both the schemes as in figure (\ref{fig:RICCA_MOVERSP_Convergence}). The convergence of RICCA is much better than that of MOVERS+.  The relative error for  RICCA drops to $1\times 10^{-12}$ uniformly.
%
%
%
%
%
\subsubsection{Horizontal slip flow}
In this test case~\cite{Manna}, a Mach $3.0$ flow slips over a Mach 2 flow with no change in pressure and density across the interface.  The computational domain considered for this test case is $\left[0, 1\right] \times \left[0, 1\right]$ with $100 \times 100$ control volumes.  The initial conditions for this test problem are as follows.
\begin{eqnarray*}
(\rho,u,v,p)_{0,y\in[0,0.5],t} =(1.4,2.0,0.0,1.0) \\
(\rho,u,v,p)_{0,y\in[0.5,1.0],t} =(1.4,3.0,0.0,1.0)
\end{eqnarray*}

\begin{figure}[h!]
    \begin{subfigure}[t]{0.3\textwidth}
        \includegraphics[width=\textwidth]{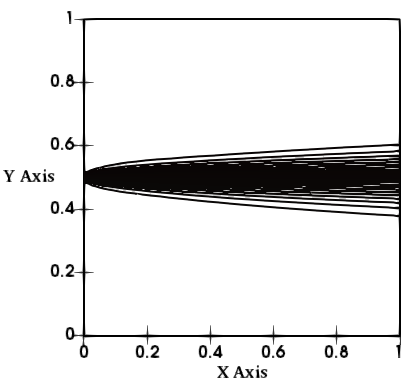}
       \caption{LLF 2O Result}
       \label{diffusedmachcontorus}
    \end{subfigure}        
    \begin{subfigure}[t]{0.3\textwidth}
        \includegraphics[width=\textwidth]{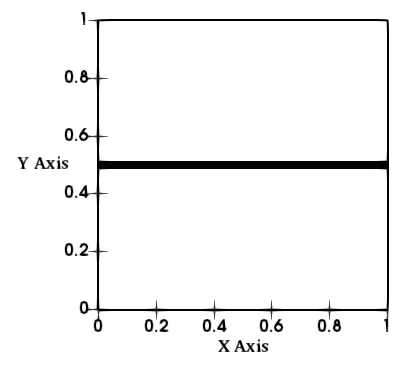}
        \caption{MOVERS+}
            \label{fig:NWSC_Slip_Flow}
    \end{subfigure}
    \begin{subfigure}[t]{0.3\textwidth}
        \includegraphics[width=\textwidth]{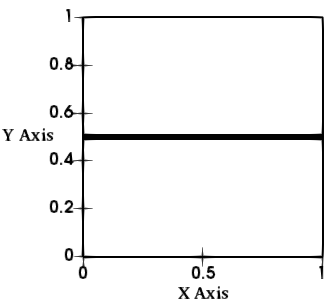}
                \caption{RICCA}
                    \label{fig:RICCA_Slip_Flow}
    \end{subfigure}
\caption{Second-order-  Mach Contours (2:0.033:3) - for Slip Flow on a $100\times100$ grid}    
\end{figure}

 Since the flow features sought are at steady state conditions, the code is run till the relative error reaches machine epsilon or the number of iterations reaches 100000.  This problem tests the accuracy of a numerical scheme in resolving a  contact discontinuity.  Many of the central and upwind schemes diffuse the contact discontinuity, due to  high numerical diffusion.  The solution of a typical diffusive scheme (Rusanov or LLF method) is shown in the figure (\ref{diffusedmachcontorus}).  Figures (\ref{fig:NWSC_Slip_Flow}) and (\ref{fig:RICCA_Slip_Flow}) show the second order accurate solutions obtained with MOVERS+ and RICCA. Both RICCA and MOVERS+ capture the grid-aligned slipstream exactly. Though not shown here, even their first order versions resolve it exactly.  

\subsubsection{Hypersonic flow past a half-cylinder}
\begin{figure}[h!]
\makebox[\textwidth][c]{%
\begin{subfigure}[b]{.3\textwidth}
\centering
\includegraphics[width=\linewidth]{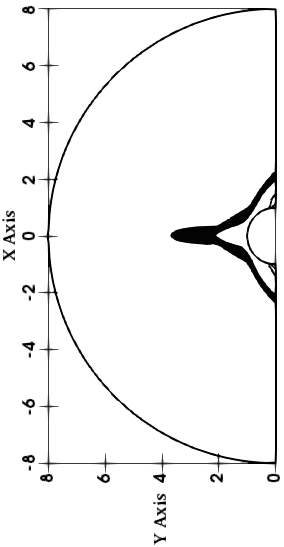}
\caption{ROE Scheme}
\label{carbuncle}
\end{subfigure}
\begin{subfigure}[b]{.3\textwidth}
\centering
\includegraphics[width=\linewidth]{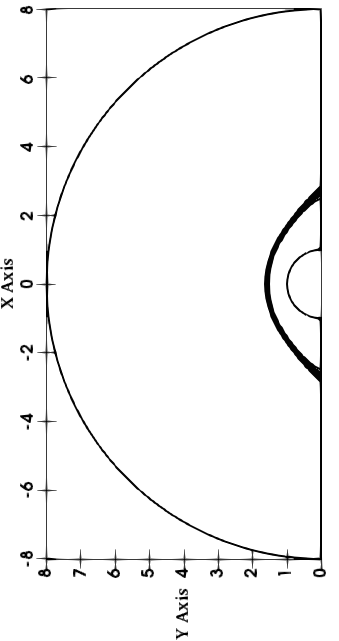}
\caption{MOVERS+}
\label{HalfCylinder_MoversP}
\end{subfigure}%
\begin{subfigure}[b]{.3\textwidth}
\centering
\includegraphics[width=0.95\linewidth]{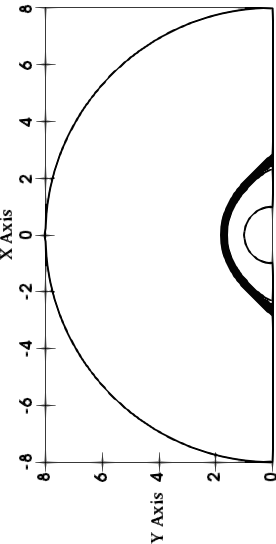}
\caption{RICCA}
\label{HalfCylinder}
\end{subfigure}
}
\caption{Comparison of 2O solutions of ROE scheme with RICCA and MOVERS+ for carbuncle effect on $240 \times 80$ Grid}
\label{fig: Compare_ROE_RICCA_SO_half_cyl}
\end{figure}

This test case assesses a numerical scheme for the shock instability called carbuncle shock discussed by Quirk~\cite{Quirk} and Meng-Sing Liou~\cite{Liou}.  A hypersonic flow of $M=6$ over a half cylinder leads to a detached bow shock in front of the bluff body.  Many of the approximate Riemann solvers like Roe scheme and low diffusive schemes produce carbuncle shocks as shown in figure (\ref{carbuncle}).  Second order (2O) accurate solutions of RICCA and MOVERS+ are compared with that of 2O solution of Roe Scheme on a $240 \times 80$ grid in figures (\ref{HalfCylinder_MoversP}) and (\ref{HalfCylinder}).  It can be observed from figure (\ref{fig: Compare_ROE_RICCA_SO_half_cyl}), that MOVERS+ and RICCA capture the bow shock without producing carbuncle shocks. MOVERS+ is shown to produce a crisper shock compared to RICCA. 

\subsubsection{Supersonic flow over forward-facing step}
\begin{figure}[htb!]
\makebox[\textwidth][c]{%
\begin{subfigure}[b]{.6\textwidth}
\centering
\includegraphics[width=\textwidth]{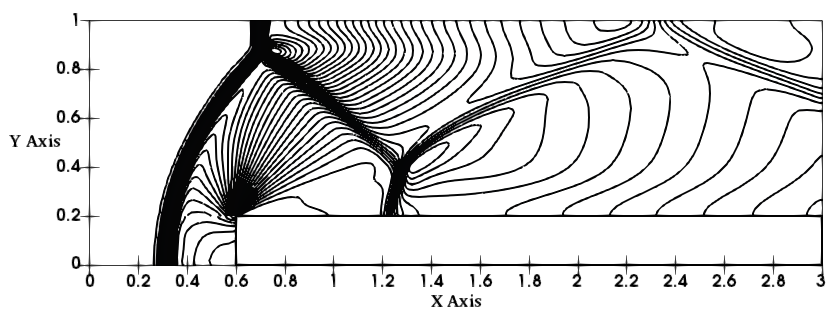}
\caption{MOVERS+ 2O accurate}
\end{subfigure}%
\begin{subfigure}[b]{.6\textwidth}
\centering
\includegraphics[width=\textwidth]{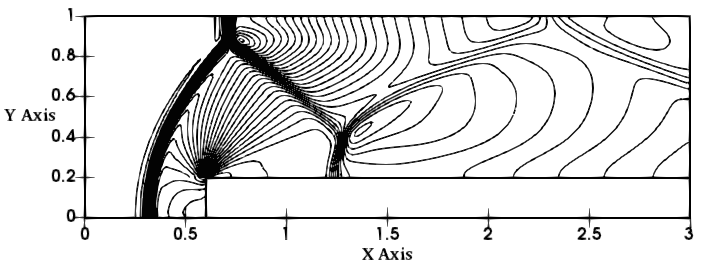}
\caption{RICCA 2O accurate}
\end{subfigure}
}
\caption{Second-order results with RICCA and MOVERS+ - Density contours (1.0:0.15:6.5) at time=4.0 - for Forward-facing step on a $240\times80$ grid}
\label{fig:Forward_facing_step_SO_240_80_ch3}
\end{figure}

In this unsteady test case~\cite{Woodward}, a Mach 3 flow enters a wind tunnel containing a forward-facing step. The computational domain ($[3,0] \times [0,1]$) consists of a step size of $0.2$ units begining at $x =0.6$. Grid sizes of  $120 \times 40$ and $240 \times 80$ are used in the simulations.  As the test case is an unsteady one, the results at time t=4.0 are presented for the fine grid.  At $t =4.0$, a lambda shock is developed and a slipstream can be seen emanating from the triple point. Results are presented with second-order accuracy in figure (\ref{fig:Forward_facing_step_SO_240_80_ch3}). Both the numerical schemes MOVERS+ and RICCA resolve the lambda shock and the reflected shocks reasonably well while MOVERS+ resolves the slipstream more accurately.  

\subsubsection{Odd-even decoupling}
\begin{figure}[h!]  
\centering
\centering
  \includegraphics[width=0.5\textwidth]{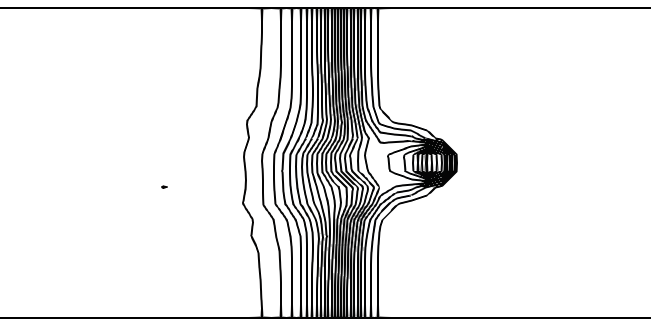}
  \caption{Odd even decoupling by 1O Roe scheme}
\label{odd_even} 
\end{figure}

\begin{figure}[h!]
\begin{subfigure}[t]{0.4 \textwidth}
\centering
  \includegraphics[width=1.23\textwidth]{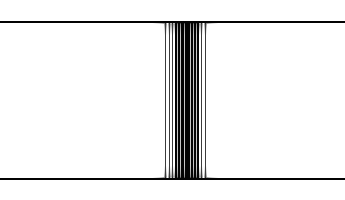}
\caption{MOVERS+ 2O Results}
  \label{oed_movers+}
\end{subfigure}
\hfill
\begin{subfigure}[t]{0.4 \textwidth}
\centering
  \includegraphics[width=\textwidth]{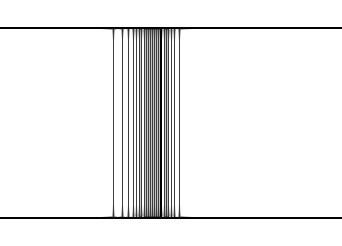}
\caption{RICCA 2O Results}
  \label{oed_ricca}
\end{subfigure}
\caption{Results for odd-even decoupling testcase, density contours on a 800$\times$20 grid at time t=100}
\label{fig: Comparison_odd_even_decoupling_ch3}
\end{figure}

This is a testcase described in ~\cite{Quirk} which assesses a numerical scheme for shock instability called odd-even decoupling. In this test case a slowly moving planar shock with Mach number $M =6$ simply travels along a long rectangular duct. The ability of the numerical schemes to avoid a numerical instability if the grid is perturbed is tested in this problem. For numerical solution, the duct is set up with a grid size of $800\times20$ unit square cells and the centerline of the grid is perturbed in the following manner:
\begin{equation*}
{y}_{i,{j}_{mid}} =
\begin{cases}
{y}_{i,{j}_{mid}} + 10^{-3} \ \text{for \emph{i} even}, \\
{y}_{i,{j}_{mid}} - 10^{-3} \ \text{for \emph{i} odd}
\end{cases}
\end{equation*}
Most of the low diffusion schemes distort the shock structure as shown in figure (\ref{odd_even}) because of the perturbed grid. For schemes like Godunov's exact Riemann solver and approximate Riemann solver of Roe, this perturbation promotes odd-even decoupling thereby destroying the planar shock structure \cite{Quirk}.   Figures(\ref{oed_movers+}) and (\ref{oed_ricca}) represent the 2O solution using MOVERS+ and RICCA.  It can be observed that the shock captured (after a long time t=100 
) using both the schemes are stable to the perturbation and no distortion of the shock structure is seen in these two schemes.

\subsubsection{Double-mach reflection (DMR)}
\begin{figure}[h!]
\begin{subfigure}[t]{0.35\textwidth}
        \includegraphics[width=\textwidth]{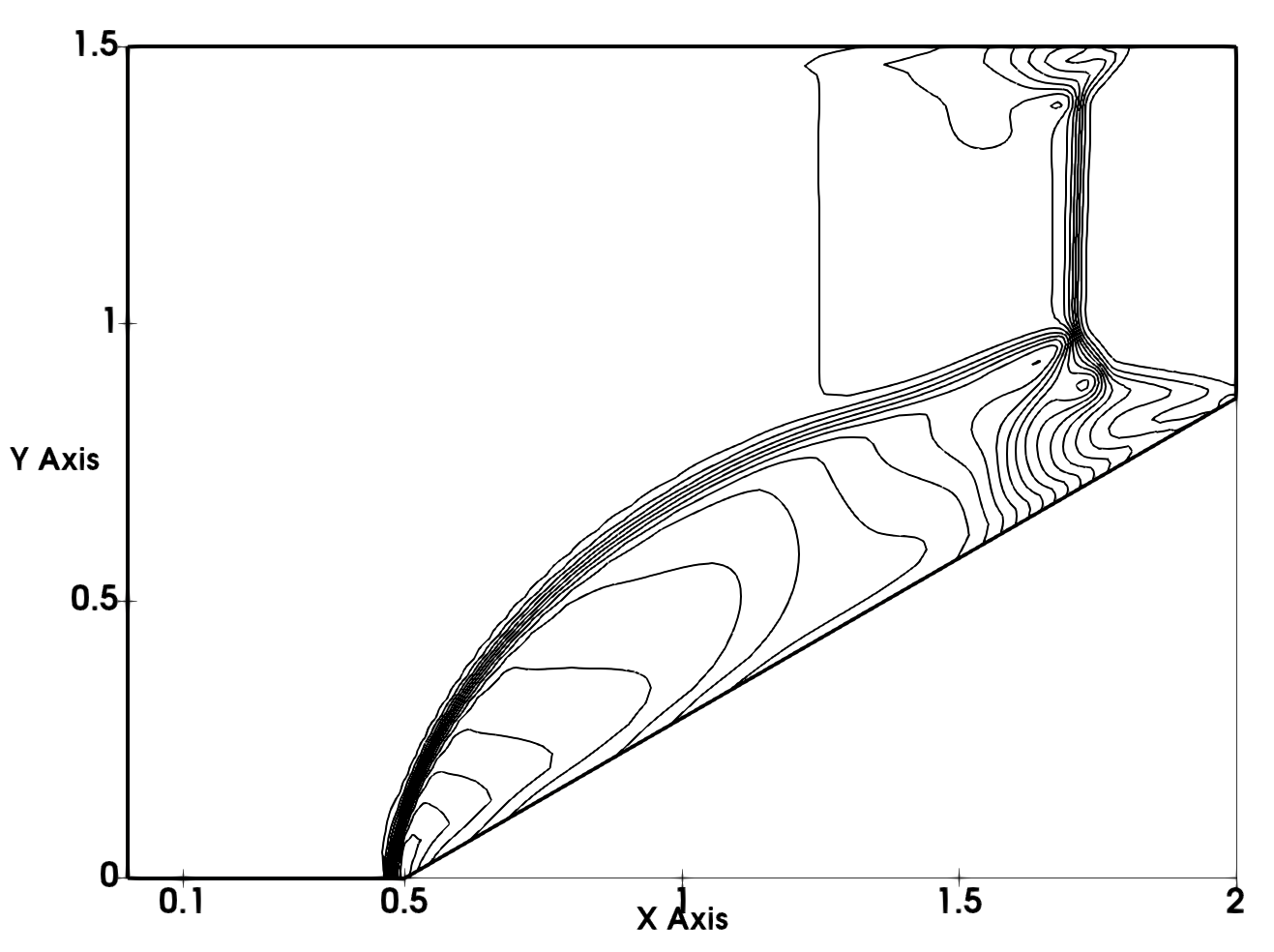}
       \caption{Kinked mach stem using 1O Roe scheme}
       \label{kinked}
\end{subfigure}
\begin{subfigure}[t]{0.35\textwidth}
\includegraphics[width=\textwidth]{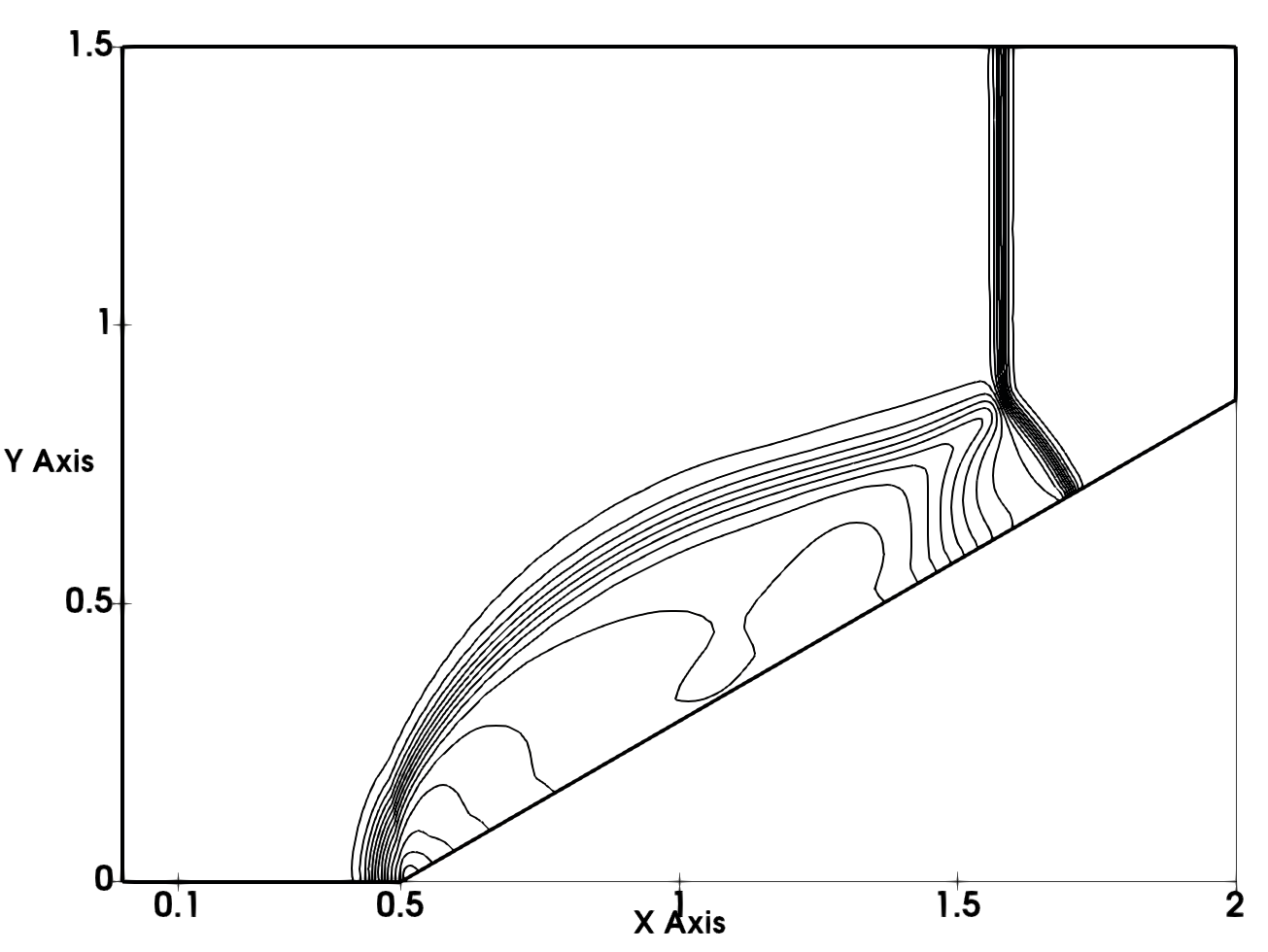}
\caption{RICCA 2O}
\end{subfigure}%
\begin{subfigure}[t]{0.35\textwidth}
\centering
\includegraphics[width=\textwidth]{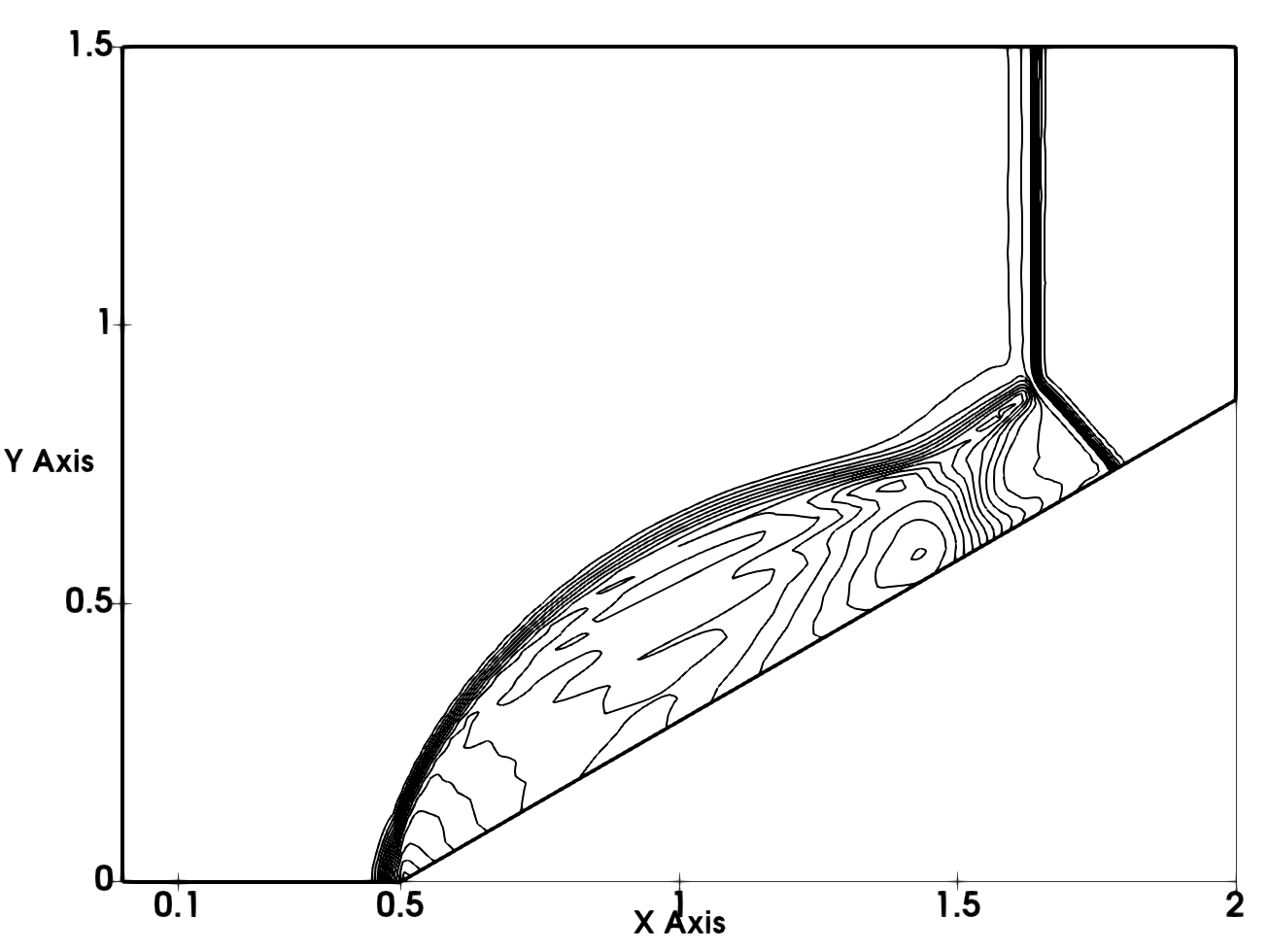}
\caption{MOVERS+ 2O}
\end{subfigure}
\caption{Density contours (0.98:0.6328:16.8) at time=0.3 - for double-mach reflection on a $240\times60$ grid
}
\label{fig:FO_SO_DMR_TestCases_ch3} 
\end{figure}

In this unsteady test case~\cite{Woodward}, a Mach 10 shock is driven down a channel
 containing a $30^\circ$ wedge. At first the simple planar shock meets the walls of the tube at right angles, but on encountering the sloping surface of the wedge, a complicated shock reflection occurs resulting in the formation of reflected shocks, Mach stems, triple points and slipstreams. A shock-instability termed kinked Mach stem~\cite{Quirk} is produced by some schemes as shown in figure (\ref{kinked}).  Results for this unsteady test case (at time t=0.3) are presented with first-order accuracy and second-order accuracy in figure (\ref{fig:FO_SO_DMR_TestCases_ch3}). Both MOVERS+ and RICCA do not produce kinked Mach stems in capturing  various features of double Mach reflection.   
%
%
%
%
%
\subsubsection{Shock Diffraction}
\begin{figure}[h!]
\begin{subfigure}[t]{0.5\textwidth}
\centering
\includegraphics[width=\linewidth]{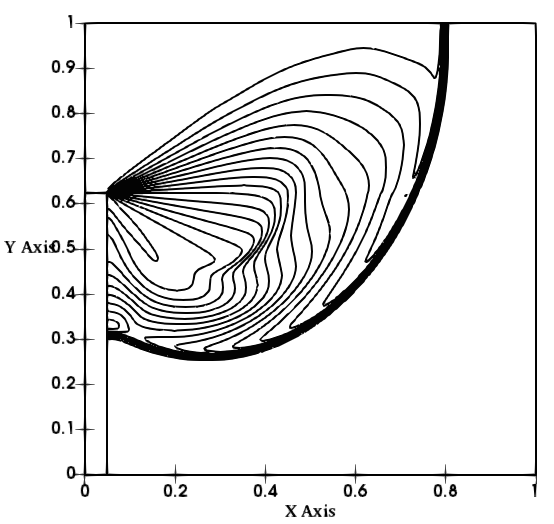}
\caption{RICCA 2O}
\label{fig: RICCA_Comparison_FO_SO_Shock_diffr_TestCases_ch3}
\end{subfigure}%
\begin{subfigure}[t]{0.5\textwidth}
\centering
\includegraphics[width=1.05\linewidth]{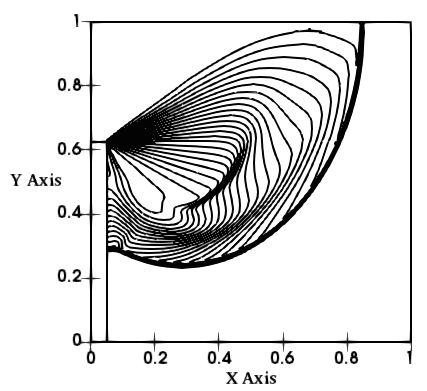}
\caption{MOVERS+ 2O}
\label{fig:NWSC_Comparison_FO_SO_Shock_diffr_TestCases_ch3}
\end{subfigure}
\caption{Density contours (0.5:0.25:6.75) at time=0.1561 - for shock diffraction test case on a $400\times400$ grid}
\label{fig:FO_SO_Shock_diffr_TestCases_ch3}
\end{figure}

This is another test case~\cite{Huang} which assesses a numerical scheme for shock instability resulting in shock anomalies and expansion shock as described in \cite{Quirk}. This test case has complex flow features involving a planar shock wave moving with incident Mach number $5.09$ 
, a diffracted shock wave around the $90^\circ$ corner and a strong expansion wave. The strong shock wave accelerates the flow and interacts with strong expansion to further complicate the flow. Other distinct flow features are a slipstream and a contact surface. Godunov and Roe schemes are known to fail for this test case~\cite{Quirk} as they admit expansion shocks without a proper fix.  Second order results for this unsteady test case (at time t=0.1561)  are presented in figures (\ref{fig: RICCA_Comparison_FO_SO_Shock_diffr_TestCases_ch3}) and (\ref{fig:NWSC_Comparison_FO_SO_Shock_diffr_TestCases_ch3}).  Both MOVERS+ and RICCA do not produce expansion shocks for this test case.  The slipstream is resolved better by  MOVERS+. 
\section{Benchmark test cases for viscous flows} \label{sec:viscous_results}
The discretization of viscous terms is as discussed in brief in section (\ref{sec:geqns}).  More details are available in \cite{Ramesh_Thesis}.      
The following 2D benchmark test cases are chosen for viscous flows such that both classical non-linear waves along with the boundary layer effects are present in evaluating the novel algorithms discussed in sections (\ref{sec:RICCA}) and  (\ref{sec:movers+}).  
\begin{itemize}
\item Flow in a viscous shock tube.
\item Interaction of oblique shockwave with laminar boundary layer.
\item Supersonic viscous flow over a thick cylindrical bump.
\end{itemize} 
\subsection{Viscous shock tube test case}
\begin{figure}[h]
\makebox[\textwidth][c]{%
\begin{subfigure}[b]{.5\textwidth}
\centering
\begin{adjustbox}{trim=0 0 0 4.6cm}%
\begin{tikzpicture}[scale = 2]
\draw [blue,very thick] (0,0) -- (4,0) --(4,1.5)--(0,1.5)--(0,0);
\draw [fill=orange] (0,0) rectangle (2,1.5);
\draw [blue,very thick] (2,0) -- (2,1.5);
\draw [very thick](1.0,0.75)node [above] {$p_4,u_4,\rho_4,T_4$} (1.5,0.75);
\draw [very thick](1.0,0.75)node [below] {Driver Section} (1.5,0.75);
\draw [very thick](1.0,1.0)node [above] {High Pressure} (1.5,1.0);
\draw [fill=gray] (2,0) rectangle (4,1.5);
\draw [very thick](3.0,0.75)node [above] {$p_1,u_1,\rho_1,T_1$} (3.5,0.75);
\draw [very thick](3.0,0.75)node [below] {Driven Section} (3.5,0.75);
\draw [very thick](3.0,1.0)node [above] {Low Pressure} (3.5,1.0);
\draw [very thick](2.0,1.7)node [above] {Diapraghm} (2.0,1.8) ;
\draw [->,red,very thick](2.0,1.7) -- (2.0,1.5) node {} ;
\end{tikzpicture}  
\end{adjustbox}
\caption{Initial Conditions in Shock tube $t=0$}
\label{shockt0}
\end{subfigure}%
\hspace{0.5cm}
\begin{subfigure}[b]{.5\textwidth}
\centering
\begin{tikzpicture}[scale = 1]
\filldraw[fill=orange, draw=black] (1,0) rectangle (2.5,3); 
\draw [very thick](1.75,1.75)node [above] {$p_4$} (1.75,1.75);
\draw [very thick](1.75,1.5)node [below] {$\rho_4$} (1.75,1.5);
\filldraw[fill=orange!70!white, draw=black] (2.5,0) rectangle (4,3); 
\draw [very thick](3.2,1.75)node [above] {$\frac{\partial p}{\partial x}<0$} (3.5,1.75);
\draw [very thick](3.2,1.05)node [above] {$\frac{\partial \rho}{\partial x}<0$} (3.5,1.05);
\filldraw[fill=orange!50!white, draw=black] (4,0) rectangle (6,3); 
\draw [very thick](4.95,1.75)node [above] {$p_3$} (4.5,1.75);
\draw [very thick](4.95,1.65)node [below] {$\rho_3$} (4.5,1.65);
\filldraw[fill=orange!30!white, draw=black] (6,0) rectangle (7.5,3); 
\draw[fill=gray, draw=black] (6.0,3.0) -- (6.0,3.9) node [above]{contact-discontinuity} ; 
\draw [very thick](6.5,1.75)node [above] {$p_2$} (6.75,1.75);
\draw [very thick](6.5,1.75)node [below] {$\rho_2$} (6.75,1.75);
\filldraw[fill=gray, draw=black] (7.5,0) rectangle (9,3); 
\draw[fill=gray, draw=black] (7.5,3.0) -- (7.5,3.35) node [above]{Shock wave} ; 
\draw [very thick](8.2,1.75)node [above] {$p_1$} (8.5,1.75);
\draw [very thick](8.2,1.65)node [below] {$\rho_1$} (8.5,1.65);
\draw[pattern=north west lines, pattern color=blue] (2.5,0) .. controls (4.0,0.25) and (5,0.25) .. (7.5,0);
\draw[pattern=north west lines, pattern color=blue] (2.5,3.0) .. controls (4.0,2.75) and (5,2.75) .. (7.5,3);
\draw [blue,very thick](4.5,0.0) -- (4.5,3.3) node [above] {Diapraghm};
\draw [ultra thick, green, latex-latex'] (2.5,3)node [above] {Rarefraction wave} -- (4.0,3.0) ;
\draw [decorate,decoration={brace,amplitude=6pt,mirror,raise=0pt},xshift=0pt,yshift=0pt](2.5,0.0) -- (3.5,0.0)node [black,midway,xshift=-4pt,yshift=-9pt] {\footnotesize LBL};
\draw [decorate,decoration={brace,amplitude=5pt,mirror,raise=0pt},xshift=0pt,yshift=0pt](4.0,0.0) -- (7.5,0.0)node [black,midway,xshift=-4pt,yshift=-09pt] {\footnotesize $p_2 = p_3$};
\draw [decorate,decoration={brace,amplitude=6pt,mirror,raise=0pt},xshift=0pt,yshift=0pt](6.5,0.0) -- (7.5,0.0)node [black,midway,xshift=0pt,yshift=-9pt] {\footnotesize LBL};
\end{tikzpicture}
\caption{Flow conditions in shock tube $t>0$}
\label{shocktg0}
\end{subfigure}%
}
\caption{Schematic of viscous flow structure inside the shock tube}
\label{shocktube_viscouseffects}
\end{figure}
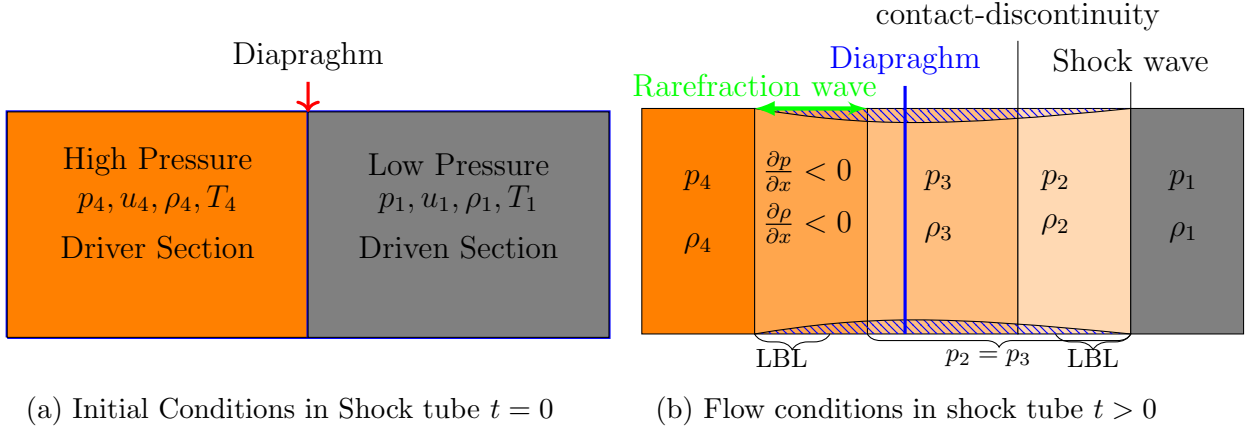

  A shock tube is a simple constant area duct with a high-pressure driver section separated from a low pressure driven section by using a diaphragm as shown in figure (\ref{shockt0}).
A comprehensive review of shock tube is presented in \cite{Olaf} and for detailed history of shock tube one can refer to \cite{krehl}.  Based on the requirement of study, the end of the driven section may be closed, open or attached to a nozzle. As the diaphragm ruptures, three waves emerge from the point of the location of the diaphragm. A shock wave gets generated and travels into the driven section, followed by a contact discontinuity and further a rarefaction wave moving in the opposite direction of the shock wave as shown in figure (\ref{shocktg0}). The strength of the shock wave depends on the pressure ratio of the gas between the driver and driven sections and its composition.  As the shock wave propagates into the driven section, all the properties of the fluid experience a discontinuous jump across a shock wave.  Across a contact-discontinuity, pressure and velocity are constant while density and temperature vary discontinuously. In the case of a rarefaction wave all the properties vary smoothly.  Apart from the propagating waves in the shock tube the viscous effects like boundary layer growth and its effects in the shock tube need attention.  The growth of an unsteady boundary layer behind the shock bringing in non-uniformity of flow across a cross-section is shown in figure(\ref{shocktg0}). This development of boundary layer accelerates the contact surface while decelerating the shock and also generates pressure waves in the duct that attenuate the shock wave, leading to a reduction in test times \cite{jameson2}. 
 Numerical simulations are carried out for viscous shock tube using MOVERS+ and RICCA. 
 
The primary objective is to check if the numerical schemes can resolve the boundary layers and the effect of these boundary layers on all the non-linear waves.  The computational domain consists of a shock tube of length $L =1$ with height to length ratio being given by $ \frac{h}{L} = 0.3$. A nonuniform grid ($141 \times 141$), stretched in the y-direction, as shown in the figure (\ref{fig:Grid_Shocktube}) is considered for the simulation.   Sod shock tube data is considered for initial conditions with Reynolds numbe, $Re =25000$ and Prandtl number, $Pr = 0.72$.

\begin{figure}[h] \label{fig:Grid_Shocktube}
\makebox[\textwidth][c]{%
\begin{subfigure}[b]{.5\textwidth}
\centering
\includegraphics[width=\textwidth,height=\textheight,keepaspectratio]{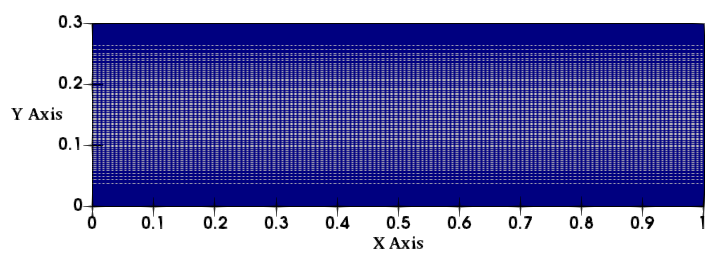}
\caption{Non-uniform grid for shock tube}
\end{subfigure}%
\begin{subfigure}[b]{.5\textwidth}
\centering
\includegraphics[scale =0.7]{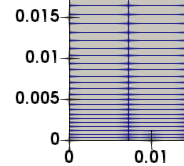}
\caption{Stretched grid near wall}
\label{fig:Grid_Shocktube}
\end{subfigure}
}
\caption{2D Grid for shock tube}
\end{figure}      
\begin{figure}[h!]
\makebox[\textwidth][c]{%
\begin{subfigure}[b]{.5\textwidth}
\centering
\includegraphics[width=\textwidth,height=\textheight,keepaspectratio]{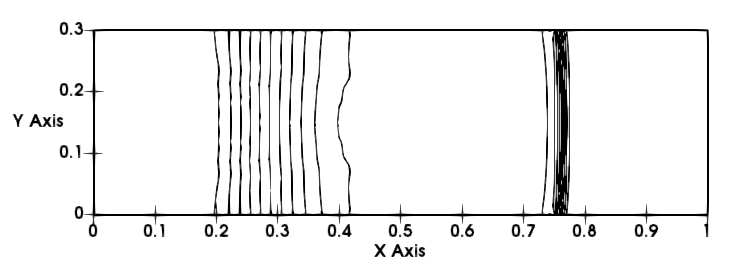}
\caption{MOVERS+}
\end{subfigure}%
\hspace{0.5mm}
\begin{subfigure}[b]{.5\textwidth}
\centering
\includegraphics[width=1.08\textwidth]{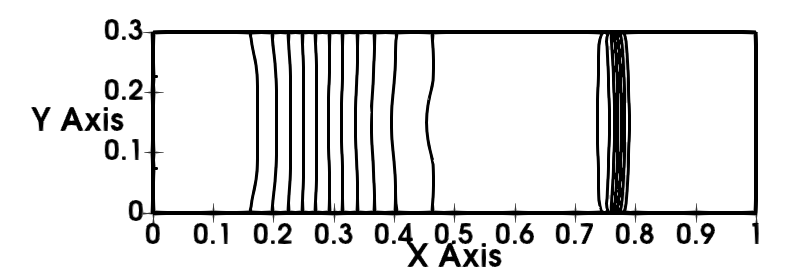}
\caption{RICCA}
\label{P_Contours}
\end{subfigure}
} 
\caption{Pressure Contours in Viscous Shock Tube} 
\label{Pressure_Contours_ShockTube}
\end{figure}

\begin{figure}[h!]
\makebox[\textwidth][c]{%
\begin{subfigure}[b]{.5\textwidth}
\centering
\includegraphics[width=1.03\textwidth,height=\textheight,keepaspectratio]{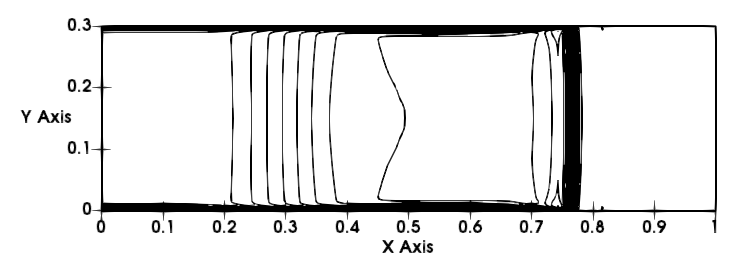}
\caption{MOVERS+}
\end{subfigure}%
\hspace{0.15mm}
\begin{subfigure}[b]{.5\textwidth}
\centering
\includegraphics[width=\textwidth,keepaspectratio]{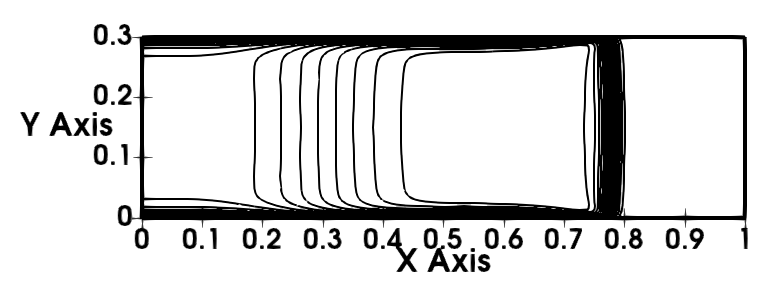}
\caption{RICCA}
\label{P_Contours}
\end{subfigure}
} 
\caption{u Velocity Contours in Viscous Shock Tube} 
\label{uvelocity_Contours_ShockTube}
\end{figure}

\begin{figure}[h!]
\makebox[\textwidth][c]{%
\begin{subfigure}[b]{.5\textwidth}
\centering
\includegraphics[width=1.03\textwidth,height=\textheight,keepaspectratio]{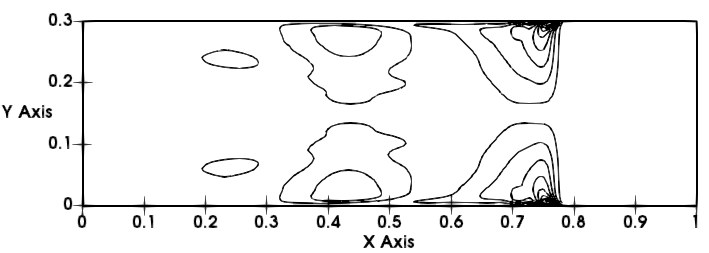}
\caption{MOVERS+}
\end{subfigure}%
\hspace{0.15mm}
\begin{subfigure}[b]{.5\textwidth}
\centering
\includegraphics[width=\textwidth,keepaspectratio]{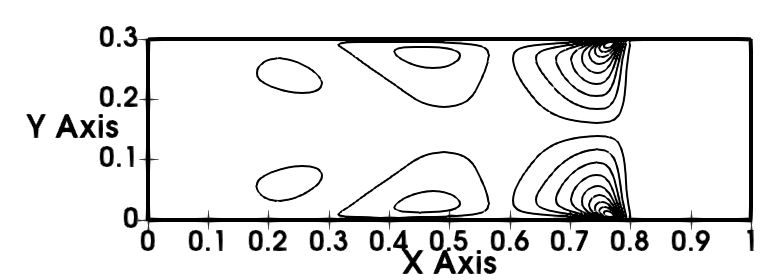}
\caption{RICCA}
\label{P_Contours}
\end{subfigure}
} 
\caption{v Velocity Contours in Viscous Shock Tube} 
\label{vvelocity_Contours_ShockTube}
\end{figure}

\begin{figure}[h!]
\makebox[\textwidth][c]{%
\begin{subfigure}[b]{.5\textwidth}
\centering
\includegraphics[width=1.03\textwidth,height=\textheight,keepaspectratio]{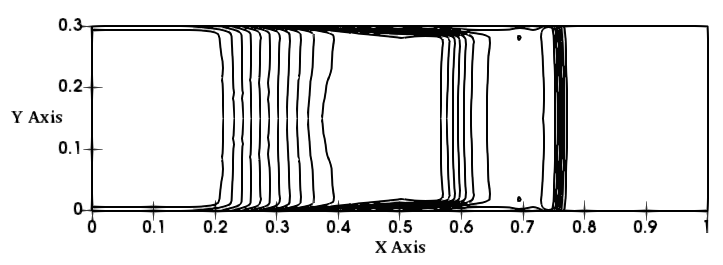}
\caption{MOVERS+}
\end{subfigure}%
\hspace{0.15mm}
\begin{subfigure}[b]{.5\textwidth}
\centering
\includegraphics[width=\textwidth,keepaspectratio]{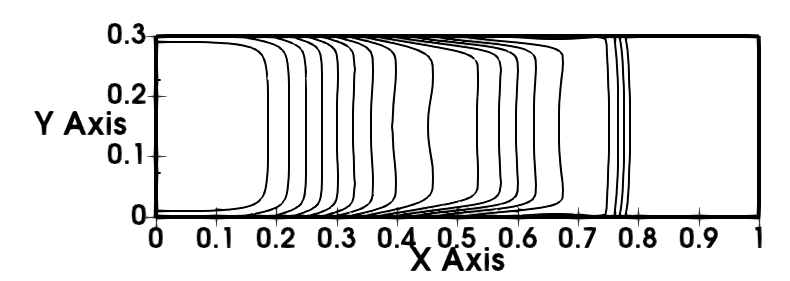}
\caption{RICCA}
\label{P_Contours}
\end{subfigure}
} 
\caption{Density Contours in Viscous Shock Tube} 
\label{rho_Contours_ShockTube}
\end{figure}

\begin{figure}[h!]
\centering
\includegraphics[scale = 0.2]{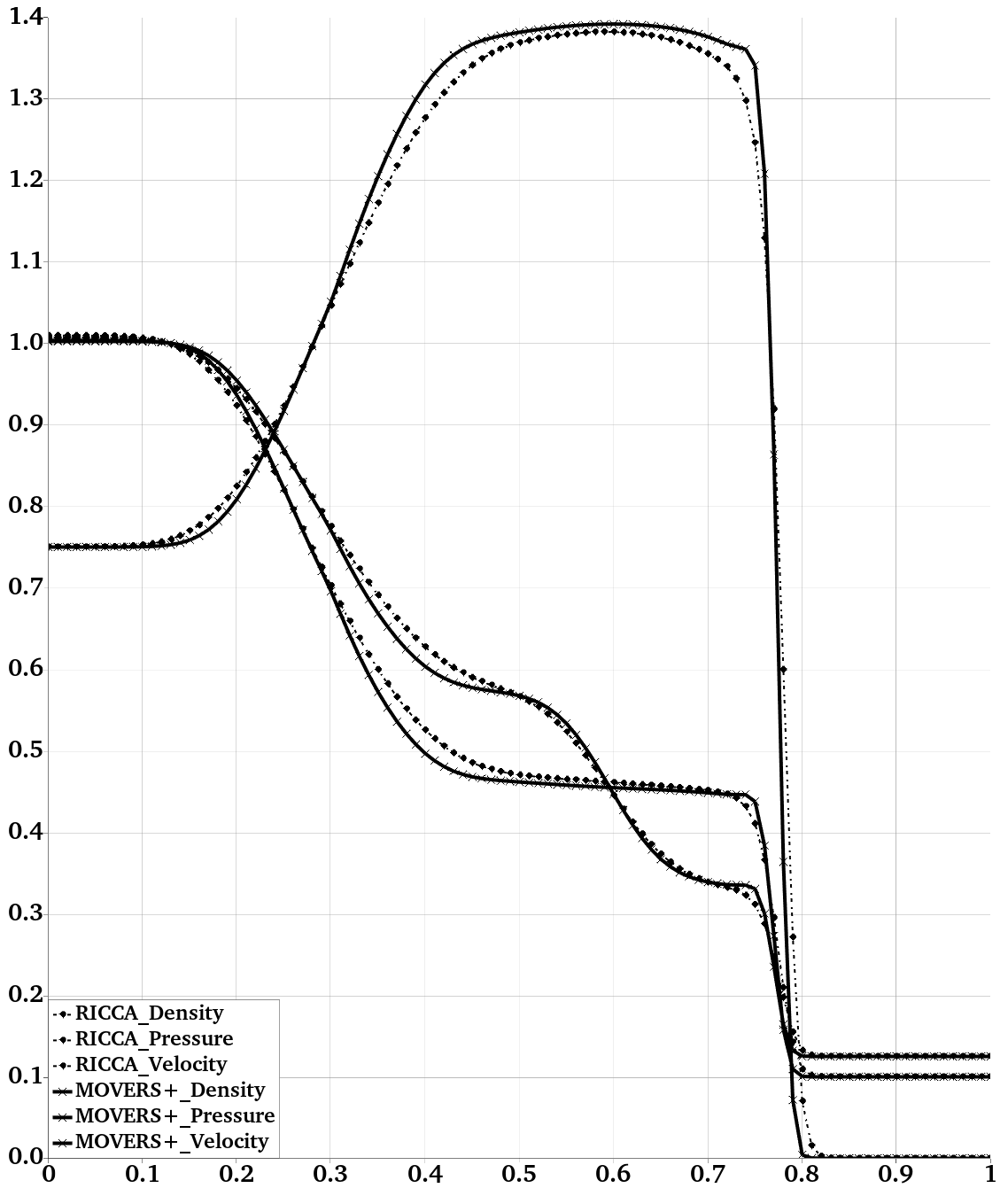}
\caption{Line Plot of Pressure, Density and Velocity Magnitude at axis using MOVERS+ and RICCA at the centreline of the shock tube}
\label{fig:CenterlineProperties}
\end{figure}


Results of second order accurate simulations with MOVERS+ and RICCA are given in figures  
 (\ref{Pressure_Contours_ShockTube}),   (\ref{uvelocity_Contours_ShockTube}) and  (\ref{fig:CenterlineProperties}).   It can observed that the shock wave and contact-discontinuity are curved and are resolved well using MOVERS+ when compared to RICCA. 
Figures (\ref{uvelocity_Contours_ShockTube},\ref{vvelocity_Contours_ShockTube}) represent the contours of u-velocity and v-velocity. Growth of boundary layer is clearly seen from these figures. Figure (\ref{fig:CenterlineProperties}) represents the variation of properties along the centre line of the shock tube. It can be observed that at the centre line the flow is behaving as an inviscid fluid and hence the variation of the properties are comparable with data of Toro test case 1. The major features that are to be observed in this test case is the development of v-velocity because of the boundary layer effects as described in \cite{jameson2}. The results presented here are not comparable with the features described in \cite{jameson2} as the schemes used by the authors in \cite{jameson2} are non-diffusive and they present results from DNS computations on an extremely fine grid. The inviscid features of the flow field are captured well and the viscous features like the curved shock, boundary layer and contours of v-velocity are resolved reasonably well. \\
\subsection{Shock wave laminar boundary layer interaction SWBLI}
\begin{figure}[h!]
  \centering
\includegraphics[width=0.5\textwidth]{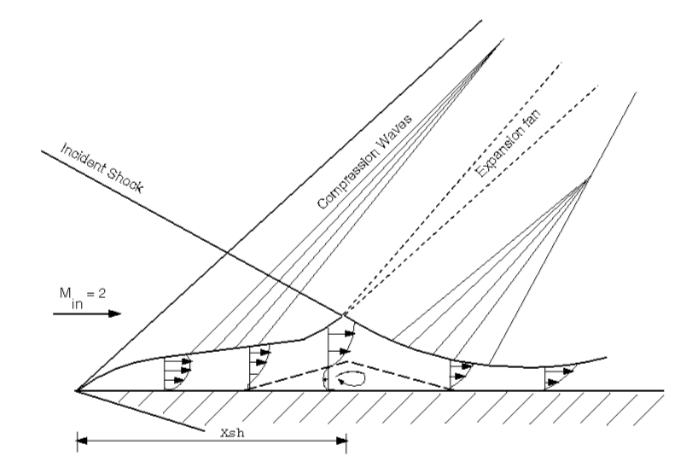}
  \caption{Schematic of oblique shockwave interacting with laminar boundary layer : Picture courtsey \cite{degrez}}
  \label{swbli}
\end{figure}

This test case represents the interaction of an oblique shock wave with a laminar boundary layer. This is a standard test case to test the ability of the numerical schemes to resolve viscous features like flow separation, bubble formation and corresponding negative skin friction. Consider a supersonic flow over a flat plate on which an oblique shock is impinged. During this process a series of events take place as depicted in the figure (\ref{swbli}). 
As the boundary layer gets developed an oblique shock evolves from the leading edge of the boundary layer and interacts with the incident oblique shock. Further as the incident oblique shock wave impinges on the laminar boundary layer, flow separation takes place inside the boundary layer. 
 Series of expansion fans and compression waves evolve from the surface of the bubble. Inside the boundary layer flow separation takes place because of flow reversal. Experimental results of such interaction was given by Hakkinen in \cite{hakkinen}, and the numerical simulations were performed by Degrez \cite{degrez}. 

The computational domain consists of region between $x \in [-0.2,1.8], y\in[0,1]$. Left boundary of the domain is initialized with supersonic inflow with Mach no $M=2.0$, Reynolds number $Re = 1\times 10^5$ and the Prandtl number $P_r = 0.72$ till $(x,y)=(0.0,0.765)$; after that the solution is initialized with corresponding post-shock conditions. The top boundary is initialized with values obtained from post-shock conditions. This ensures an oblique shock resulting at an angle of $30^{\circ}$. The bottom boundary has symmetry condition for $x \in [-0.2,0.0]$ and no-slip wall boundary conditions from $x\in[0.0,1.8]$. The right boundary is treated as a supersonic exit. 

\begin{figure}
\makebox[\textwidth][c]{%
\begin{subfigure}[b]{.5\textwidth}
\centering
\includegraphics[width=\textwidth,height=\textheight,keepaspectratio]{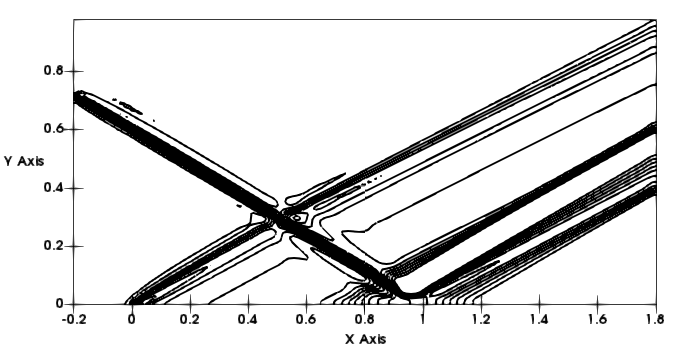}
\caption{MOVERS+}
\label{fig:SWBLI_Pressure_141_201_MOVERS_NWSC}
\end{subfigure}
\begin{subfigure}[b]{.5\textwidth}
\centering
\includegraphics[width=\textwidth]{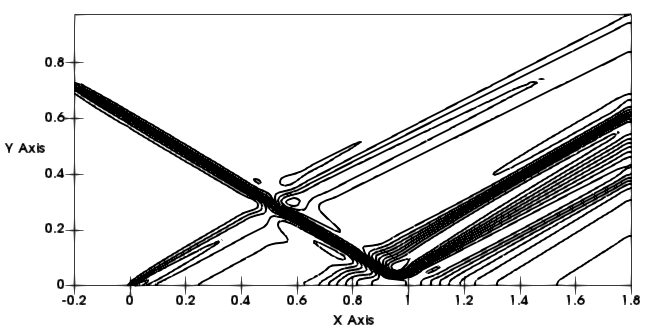}
\caption{RICCA}
\label{fig:SWBLI_Pressure_141_201_RICCA}
\end{subfigure}
}
\caption{Pressure contours}
\label{fig:PContours}
\end{figure}

\begin{figure}  
\makebox[\textwidth][c]{%
\begin{subfigure}[b]{.5\textwidth} 
\label{recirculation_zone_streamlines}
\centering
\includegraphics[width=\textwidth]{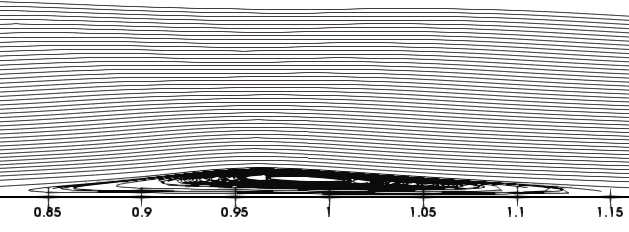}
\caption{Streamlines in recirculation zone}
\end{subfigure}
\begin{subfigure}[b]{.5\textwidth} 
\label{fig:Recirculation_Zone_Streamlines} 
\centering
\includegraphics[width=\textwidth]{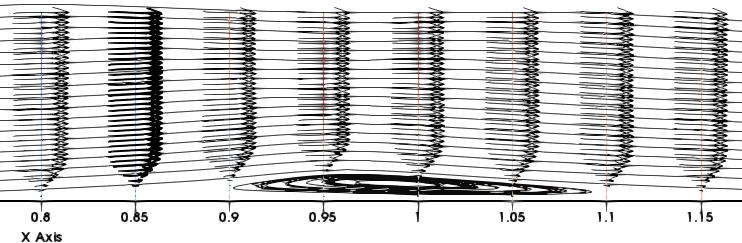}
\caption{Velocity vectors overlapped with streamlines at recirculation zone}
\end{subfigure}
}
\caption{Recirculation zone: streamlines and velocity vectors} 
\label{fig:Velocity_Vectors}
\end{figure}

The results for SWBLI using MOVERS+ and RICCA are presented here. Figure (\ref{fig:PContours}) represents the pressure contours from MOVERS+ and RICCA. It can be observed  that the incident shock on to the boundary layer, the leading edge shock from the boundary layer, the reflected shock, the expansion fans and the recirculation bubble in the boundary layer are well resolved. The resolution of shocks is good in the case of MOVERS+ when compared to the RICCA.  Further, the figure (\ref{fig:Velocity_Vectors}) shows the streamlines in the recirculation zone in the separation bubble and  velocity vectors in the recirculation zone.  

%
%
%
%


%
%
%
%
%

\begin{figure}[h!]
\makebox[\textwidth][c]{%
\begin{subfigure}[b]{.5\textwidth}
\centering
\includegraphics[width=\textwidth]{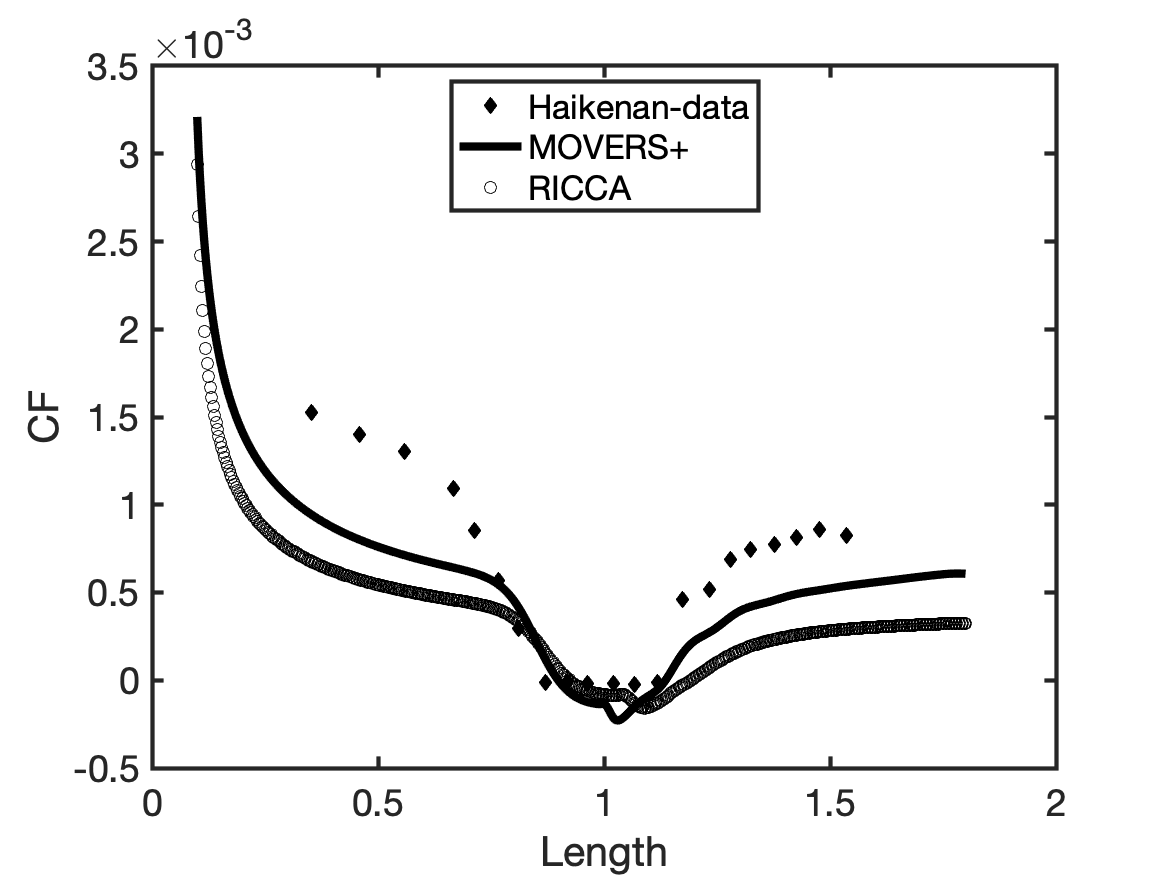}
\caption{Skin Friction Coefficient}
\label{CF_SWBLI}
\end{subfigure}
\begin{subfigure}[b]{.5\textwidth}
\centering
\includegraphics[width=\textwidth,height=0.29\textheight]{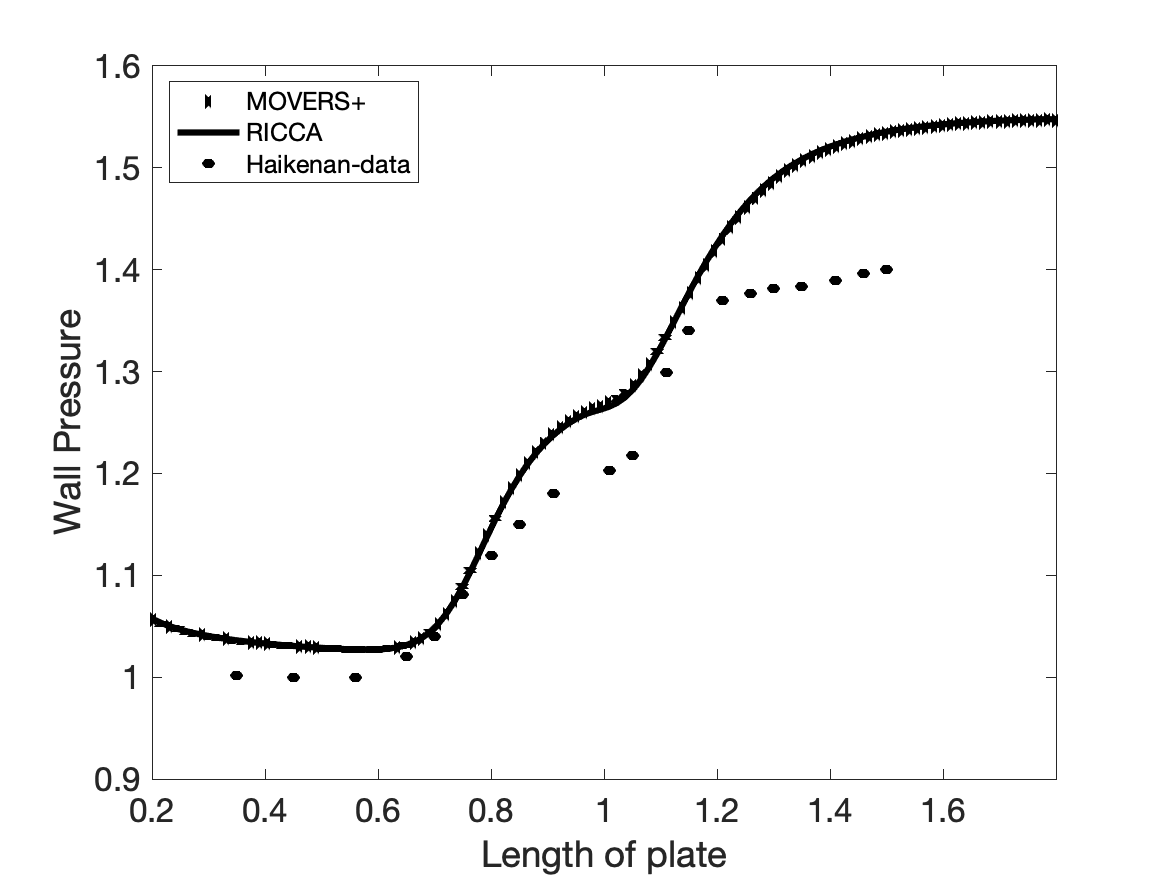}
\caption{Wall Pressure}
\label{Wall_Pressure}
\end{subfigure}
}
\caption{Comparison of coefficient of skin friction and wall pressure on wall.}
\end{figure}

The comparison of wall pressures, given by $\frac{p_w}{p_{\infty}}$, obtained with MOVERS+, RICCA and the experimental data of Hakkinen in \cite{hakkinen} is shown in figure (\ref{Wall_Pressure}).  Both MOVERS+ and RICCA follow the trendline of the wall pressure but some deviation from the experiments is found which needs further 
Figure (\ref{CF_SWBLI}) refers to the coefficient of skin friction on the wall surface. It can be observed that, while both MOVERS+ and RICCA resolve the negative skin friction, MOVERS+ gets the negative value of skin friction much better than RICCA. Resolving the separation bubble and correspondingly obtaining the negative values of skin friction is very difficult and is possible only with second order accuracy together with basically low diffusive schemes and a proper grid which resolves the boundary layer well.
The deviation of the skin friction coefficient away from the bubble region in both the schemes is probably due to the numerical diffusion present in the numerical schemes in smooth regions.  

\subsection{Viscous flow over a cylindrical bump in a channel}
The third test case considered is a low supersonic viscous flow over a $4\%$ thick circular bump in a channel as given in \cite{Parthasarathy}. As the supersonic flow approaches the bump, shock evolves from the leading edge of the bump and reflects from the top wall. This reflected shock impinges on the boundary layer of the bottom wall and a flow separation phenomenon can be observed based on the strength of the shock.
For the incoming flow, Reynolds number of $Re=8000$ and Mach no of $M=1.4$ are prescribed. The computational domain considered here is $[0,3]\times[0,1]$, where a circular bump of chord length $1$ is placed at at $[1,0]$.
 The following boundary conditions are imposed. The left boundary is considered as a supersonic inlet with $M_{\infty} = 1.4$ and the domain is initialized with the free stream conditions. The bottom part of the domain with a cylindrical bump is considered as a wall with the no-slip conditions for $x\in[1,3]$.  Symmetry conditions are enforced in the bottom boundary for  $x\in[0,1]$.  The right boundary is supersonic exit boundary and the top boundary is considered as an inviscid wall and hence the flow tangency condition is imposed.  A stretched grid size of $160 \times 160 $ is considered for numerical simulations. The grid is stretched in the y direction in order to capture viscous features  effectively.

\begin{figure}[h!]
\makebox[\textwidth][c]{%
\begin{subfigure}[b]{.5\textwidth}
\centering
\includegraphics[width=\textwidth]{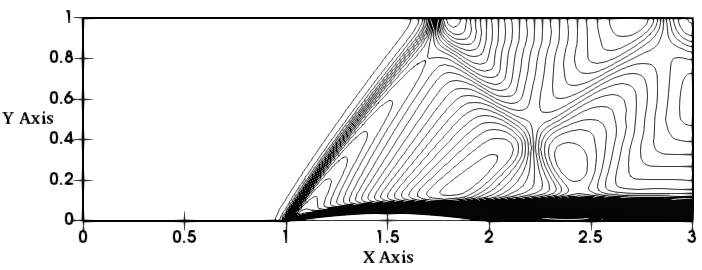}
\caption{RICCA}
\end{subfigure}
\begin{subfigure}[b]{.5\textwidth}
\centering
\includegraphics[width=\textwidth]{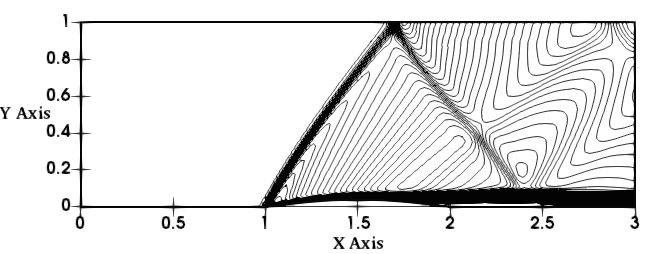}
\caption{MOVERS+}
\end{subfigure}
}
\caption{Mach Contours $[0:0.0215:1.5]$}
\label{fig:FOB_MachContorus}
\end{figure}

 Numerical simulations are carried out using MOVERS+ and RICCA. It can be seen from the figure (\ref{fig:FOB_MachContorus}) that a shock wave evolves from the leading edge of the bump. Since the top wall has only flow tangency condition the shock wave just reflects from the wall. At the bottom wall from the leading edge no-slip boundary condition gets enforced hence the boundary layer starts growing. As the flow encounters the trailing edge of the bump it separates where the reflected shock impinges on the boundary layer and gets reattached in the downstream direction. 
 These flow features are captured accurately by MOVERS+ while the numerical diffusion in RICCA leads to inaccurate resolution.  
  It can be seen that RICCA is more diffusive as is evident also in the figure (\ref{fig:FOB_MachContorus}). It can also be observed that RICCA does not show the separation region distinctively where as MOVERS+ shows the separation region.
\begin{figure}[htb!]
\centering
\includegraphics[scale =0.4]{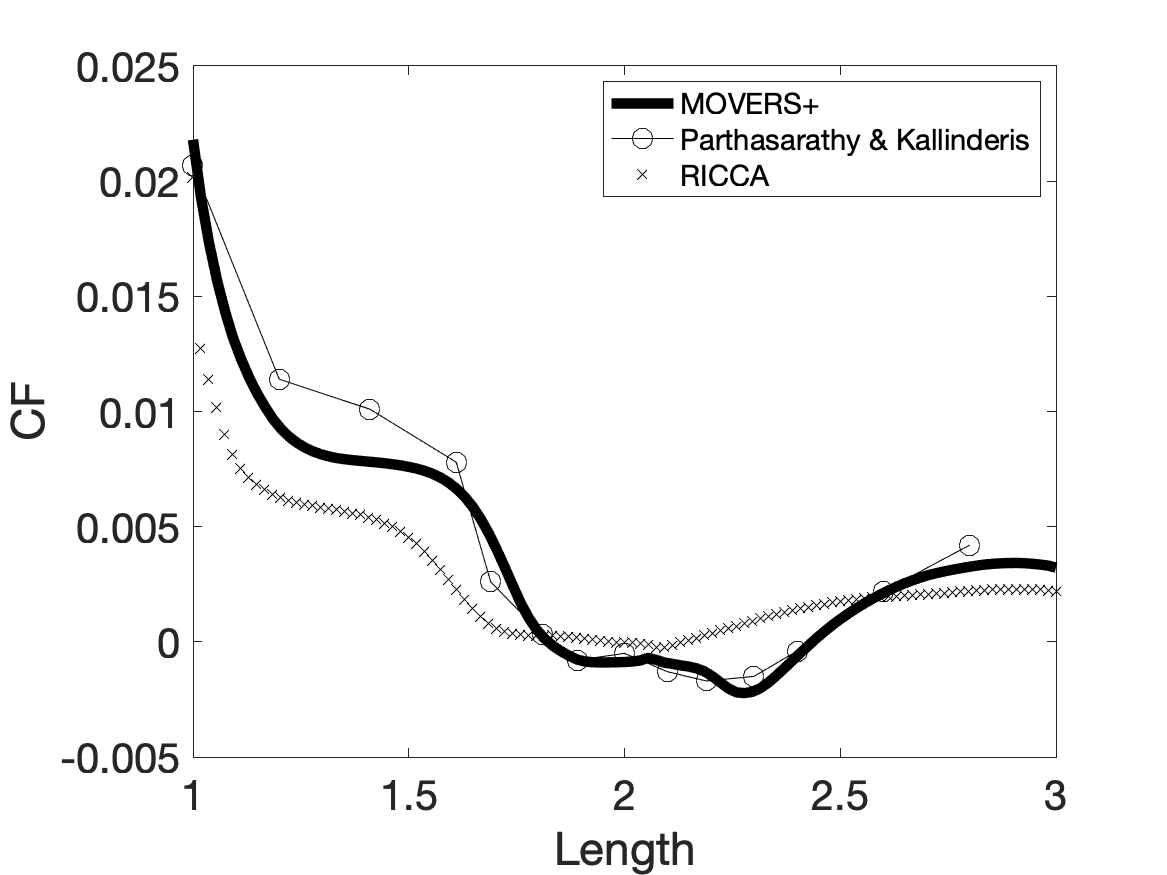}
\caption{Comparison of coefficient of friction on wall surface}
\label{fig:CF_FOB}
\end{figure}
 Figure (\ref{fig:CF_FOB}) shows a comparison of the skin friction coefficient with the data taken from \cite{Parthasarathy}. A very good match of skin friction coefficient is obtained with MOVERS+ where as skin friction plot of RICCA is deviating from the reference plot.\\

\section{Summary}
In this work, two novel algorithms are presented. The first algorithm, RICCA, is based on generalized Riemann invariants for a contact discontinuity. This leads to the coefficient of numerical diffusion equalling the fluid velocity. This diffusion helps in resolving steady contact discontinuities exactly but is insufficient in resolving the shocks in a robust way. Therefore, an additional numerical diffusion based on sound speed is added, which is scaled with the sign function of pressure difference so that it goes to zero near steady contact-discontinuity. The second algorithm is based on a modification of the previously introduced numerical diffusion matching RH conditions in a simple central solver. The modification removes the wave speed correction present in the previous version and leads to a very accurate scheme but too low in numerical diffusion. Hence an additional diffusion, taken from RICCA, is added with the help of a shock sensor. The resulting scheme, MOVERS+, is accurate and yet robust. The formulations for both the numerical schemes are simple, do not need any wave speed correction, independent of eigen-structure and do not need entropy fix. These numerical schemes work well for various benchmark test cases involving shock instabilities and shock anomalies. Both the numerical schemes are capable of resolving steady contact-discontinuities exactly. Further, the first order solutions obtained by MOVERS+ on any grid are comparable with the second order results of RICCA. MOVERS+, though does not capture steady shock exactly, has better shock capturing capabilities than RICCA. RICCA is less diffusive than the LLF scheme and is capable of capturing steady contacts exactly and hence can be a better alternative to LLF scheme. Further, numerical simulations of viscous 2D flows have been carried out using MOVERS+ and RICCA. As the boundary layer can be considered as a slipstream, these numerical schemes are expected to capture them accurately. However, MOVERS+ performed much better than RICCA for viscous flows especially in resolving flow separation bubbles. While the exact contact-discontinuity capturing feature of RICCA for inviscid flows led to accurate resolution of inviscid features, the excessive numerical diffusion in smooth regions has clearly affected the resolution of viscous regions.  

\section*{Acknowledgment} The third author thanks Prof. Francois Dubois for some interesting and fruitful discussions.  





\end{document}